\documentclass[11pt,a4paper,english]{smfbook}
\usepackage{amssymb}
\usepackage{latexsym}
\usepackage{xypic}
\usepackage{eufrak}
\usepackage{euscript}
\usepackage{amscd}

\numberwithin{equation}{chapter}
\newtheorem{thm}{Theorem}[chapter]
\newtheorem{defi}{Definition}[chapter]
\newtheorem{lem}{Lemma}[chapter]

\newtheorem{cor}{Corollary}[chapter]

\newtheorem{rema}{Remark}[chapter]
\newtheorem{notation}{Notation}[chapter]
\def\ni{\noindent}
\def\di{\displaystyle}

\font\tenbb=msbm10
\def\cal{\mathcal}
\def\stoc{\EuScript{S}}

\font\aa=ygoth
\def\bab{\hbox{\aa d}}
\def\ima{\hbox{\rm Im}}
\def\rea{\hbox{\rm Re}}
\def\cC{\hbox{\tenbb C}}
\def\rR{\hbox{\tenbb R}}
\def\nN{\hbox{\tenbb N}}

\def\grad{\hbox{\rm grad}}

\def\kk{\mathcal A}
\def\Cu{\EuScript{C}^1}
\def\Nd{\EuScript{N}^1}

\renewcommand{\d}{\displaystyle}
\renewcommand{\geq}{\geqslant}
\renewcommand{\leq}{\leqslant}

\newcommand{\N}{\mathbb{N}}

\newcommand{\R}{\mathbb{R}}

\newcommand{\C}{\mathbb{C}}

\begin{document}
\setcounter{tocdepth}{2}
\baselineskip 6mm
\title{Stochastic embedding of dynamical systems}
\author{Jacky CRESSON}
\address{Universit\'e de Franche-Comt\'e, \'Equipe de Math\'ematiques de Besan\c{c}on, CNRS-UMR 6623, 16 route de Gray, 25030
Besan\c{c}on cedex, France.}
\email{cresson@math.univ-fcomte.fr}
\author{S\'ebastien DARSES}
\address{Universit\'e de Franche-Comt\'e, \'Equipe de Math\'ematiques de Besan\c{c}on, CNRS-UMR 6623, 16 route de Gray, 25030
Besan\c{c}on cedex, France.}
\email{darses@math.univ-fcomte.fr}

\date{July 2, 1991}
\subjclass{Primary 54C40, 14E20;\\Secondary 46E25, 20C20}
\keywords{\texttt{Stochastic calculus, Dynamical systems, Lagrangian systems, Hamiltonian systems.}}

\maketitle

\begin{abstract}
Most physical systems are modelled by an ordinary or a partial differential equation, like the
$n$-body problem in celestial mechanics. In some cases, for example when studying the long term behaviour of the solar
system or for
complex systems, there exist elements which can influence the dynamics of the system which are not well modelled or
even known. One way to take these problems into account consists of looking at the dynamics of the system on a larger class
of objects, that are eventually stochastic. In this paper, we develop a theory for the stochastic embedding of ordinary
differential equations. We apply this method to Lagrangian systems. In this particular case, we extend many results of
classical mechanics namely, the least action principle, the Euler-Lagrange equations, and Noether's theorem. We also obtain
a Hamiltonian formulation for our stochastic Lagrangian systems. Many applications are discussed at the end of the paper.
\end{abstract}

\newpage
\tableofcontents

\newpage
\chapter*{Introduction}

Ordinary  as well as partial differential equations play
a fundamental role in most parts of mathematical physics. The story begins with
Newton's formulation of the law of attraction and the corresponding
equations which describe the motion of mechanical systems. Regardless the
beauty and usefulness of these theories in the study of many important natural
phenomena, one must keep in mind that they are based on experimental facts, and
as a consequence are only an approximation of the real
world. The basic example we have in mind is the motion of the planets in the solar
system which is usually modelled by the famous $n$-body problem, i.e. $n$
points of mass $m_i$ which are only submitted to their mutual gravitational
attraction. If one looks at the behaviour of the solar system for finite time
then this model is a very good one. But this is not true when one looks at the
long term behaviour, which is for instance relevant when dealing with the so called {\it chaotic}
behaviour of the solar system over billions years, or when trying to predict ice ages
over a very large range of time. Indeed, the $n$-body problem is a conservative
system (in fact a Lagrangian system) and many non-conservative effects, such as
tidal forces between planets, will be of increasing importance along the
computation. These non-conservative effects push the model outside the category of Lagrangian
systems. You can go further by considering effects due to the changing in the oblateness
of the sun. In this case, we do not even know how to model such kind of
perturbations, and one is not sure of staying in the category of differential
equations\footnote{Note that in the context of the solar system we have two different problems: first, if one
uses only Newton's gravitational law, one must take into account the entire universe to model the behaviour of the planets.
This by itself is a problem which can be studied by using the classical perturbation theory of
ordinary differential equations. This is different if we want to speak of the ``real" solar system for which we must consider
effects that we ignore. In that case, even the validation of the law of gravitation as a real law of nature is not clear.
I refer to \cite{cr4} for more details on this point.}.\\

As a first step, this paper proposes tackling this problem by introducing a
natural stochastic embedding procedure for ordinary or partial differential
equations. This consists of looking for the behaviour of stochastic processes submitted to
constraints induced by the underlying differential equation\footnote{This strategy is part of a general programme called
the embedding procedure in \cite{cr3} and which can be used to embed ordinary differential equations not only
on stochastic processes but on general functional spaces. A previous attempt was made in \cite{cr1},\cite{cr2} in the
context of the non-differentiable embedding of ordinary differential equations.}. We point out that this strategy is
different from the standard
approach based on stochastic differential equations or stochastic dynamical systems,
where one gives a meaning to ordinary differential equations perturbed by a
small random term. In our work, no perturbations of the underlying equation are carried out.\\

A point of view that bears some resemblance to ours is contained in V.I. Arnold's {\it materialization of
resonances} (\cite{ar2},p.303-304), whose main underlying idea can be briefly explained as follows: the divergence of the
Taylor expansion of the $\arctan x$ function at $0$ for $\mid x \mid >1$ can be proved by computing the
coefficients of this series. However, this does not explain the reason for this divergence behaviour. One can obtain a better
understanding by extending the function to the complex plane and by looking at its singularities at $\pm i$.
The same idea can be applied in the context of dynamical systems. In this case, we look for the obstruction to
linearization of a real systems in the complex plane. Arnold has conjectured that this is due to the
accumulation of periodic orbits in the complex plane along the real axis. In our
case, one can try to understand some properties of the trajectories of
dynamical systems by using a suitable extension of its domain of definition. In our work, we
give a precise sense to the concept of differential and partial differential equations in the class of
stochastic processes. This procedure can be viewed as a first step toward the general ``stochastic programme" as described
by Mumford in \cite{mum}.\\

Our embedding procedure is based on a simple idea: in order to write down differential or
partial differential equations, one uses derivatives. An ordinary differential
equation is nothing else but a differential operator of order one\footnote{In this
case, we can also speak of vector fields.}. In order to embed ordinary
differential equations, one must first extend the notion of derivative so that it makes sense
in the context of stochastic processes. By extension, we mean that
our stochastic derivative reduces to the classical derivative for deterministic
differentiable processes. Having this extension, one easily defines in a unique
way, the stochastic analogue of a differential operator, and as a consequence,
a natural embedding of an ordinary differential equation on
stochastic processes.\\

Of course, one can think that such a simple procedure will not produce anything
new for the study of classical differential equations. This is not the case.
The main problem that we study in this paper is the embedding of natural
Lagrangian systems which are of particular interest for classical mechanics. In this context, we obtain some
numerous surprising results, from the existence of a coherent least action principle with respect to the stochastic embedding
procedure, to a derivation of a stochastic Noether theorem, and passing by a new derivation of the
Schr\"odinger equation. All these points will be described with details in the following.\\

Two companion papers (\cite{ccbd1},\cite{ccbd2}) give an application of this
method to derive new results on the formation of planets in a protoplanetary
nebulae, in particular a proof of the existence of a so called Titus-Bode law
for the
spacing of planets around a given star.
\vskip 7mm
The plane of the paper is as follow:\\

In a first part, we develop our notion of a stochastic derivative and study in details all its properties.\\

Chapter \ref{chapnelson} gives a review of the stochastic calculus developed by Nelson \cite{ne1}. In particular, we
discuss the classical definition of the backward and forward Nelson derivatives, denoted by $D$ and $D_*$, with respect to
dynamical problems. We also define a class of stochastic process called good diffusion processes for which one can
compute explicitly the Nelson derivatives.\\

In Chapter \ref{chapstocder} we define what we call an abstract extension of the classical derivative. Using the
Nelson derivatives, we define an extension of the ordinary derivative on stochastic processes, which we call the
{\it stochastic derivative}. As pointed out previously, one imposes that the stochastic derivative reduces to the classical
derivative on differentiable deterministic processes. This constraint ensures that the
stochastic analogue of a PDE  contains the classical PDE. Of course such a
gluing constraint is not sufficient to define a rigid notion of stochastic
derivative. We study several natural constraints which allow us to obtain a unique extension of the classical derivative
on stochastic processes as
\begin{equation}
{\cal D}_{\mu} =\di {D+D_* \over 2} +i\mu \di {D-D_* \over 2}, \ \mu =\pm 1 .
\end{equation}
By extending this operator to complex valued stochastic processes, we are able to define the iterate of ${\cal D}$, i.e.
${\cal D}^2 ={\cal D}\circ {\cal D}$ and so on. The main surprise is that the real part of ${\cal D}^2$ correspond to
the choice of Nelson for acceleration in his dynamical theory of Brownian motion. However, this result depends on
the way we extend the stochastic derivative to complex valued stochastic processes. We discuss several alternative which
covers well known variations on the Nelson acceleration.\\

In Chapter \ref{chapprop} we study the product rule satisfied by the stochastic derivative which is a fundamental ingredient
of our stochastic calculus of variation. We also introduce an important class of stochastic processes, called Nelson
differentiable, which have the property to have a real valued stochastic derivative. These processes play a fundamental
role in the stochastic calculus of variation as they define the natural space of variations for stochastic processes.
\vskip 6mm
The second part of this article deals specifically with the definition of a stochastic embedding procedure for ordinary
differential equations.\\

Chapter \ref{chapembed} associate to a differential operator of a {\it given form} acting on sufficiently regular
functions a unique operator acting on stochastic processes and defined simply by replacing the classical derivative by
the stochastic derivative. This is this procedure that we call the stochastic embedding procedure. Note that the form of
this procedure acts on differential operators of a given form. Although the procedure is canonical for a given form
of operator, it is not canonical for a given operator.\\

The previous embedding is formal and does not take constraints which are of dynamical nature, like the reversibility of
the underlying differential equation. As reversibility plays a central role in physics, especially in celestial mechanics
which is one domain of application of our theory, we discuss this point in details. We introduce an embedding which respect
the reversibility of the underlying equation. Doing this, we see that we must restrict attention to the real part of our
operator, which is the unique one to possess this property in our setting. We then recover under dynamical and algebraic
arguments studies dealing with particular choice of stochastic derivatives in order to derive quantum mechanics from classical
mechanics under Nelson approach.
\vskip 6mm
The third part is mainly concerned with the application of the stochastic embedding to Lagrangian systems.\\

We consider autonomous\footnote{This restriction is due to technical difficulties.} Lagrangian systems $L(x,v)$,
$(x,v)\in U\subset \R^d \times \R^d$, where $U$ is an open set, which satisfy a number of conditions, one of it being that it must be holomorphic with respect to the second variable which
represent the derivative of a given function. Such kind of Lagrangian functions are called admissible. Using the
stochastic embedding procedure we can associate to the classical Euler-Lagrange equation a stochastic one which has the
form
$$
\di {\partial L\over \partial x} (X(t),{\cal D} X (t))=\di {\cal D}
\left [ \di
{\partial L\over \partial v} (X(t),{\cal D}X(t))
\right ] ,
\eqno{(SEL)}
$$
where $X$ is a real valued stochastic process.\\

At this point, our manipulation is only formal and one can ask if this embedding is significant or not. We then remark that
the Lagrangian function $L$ keep sense on stochastic processes and can be considered as a functional. As a consequence,
we can search for the existence of a least action principle which gives the stochastic Euler-Lagrange equation (SEL). The
existence of such a stochastic least action principle is far from being trivial with respect to the embedding procedure.
Indeed, it must follows from a stochastic calculus of variations which is not developed apart from this procedure. Our
problem can then be formalize as the following diagram:
\begin{equation}
\begin{CD}
L (x,dx/dt ) @> \mbox{\rm LAP} >> EL\\
@VV{{\cal S}}V @VV{{\cal S}}V\\
L(X,{\cal D} X ) @> \mbox{\rm SLAP}\ ? >> (SEL) ,
\end{CD}
\end{equation}
where $LAP$ is the least action principle, $\cal S$ is the stochastic embedding procedure, (EL) is the classical Euler-Lagrange equation associated to $L$ and
$SLAP$ the at this moment unknown stochastic least action principle. The existence of such a principle is called the
coherence problem.\\

Chapter \ref{chapvar} develop a stochastic calculus of variations for functionals of the form
\begin{equation}
E\left [
\int_a^b L(X(t),{\cal D}X(t))\, dt \right ] ,
\end{equation}
where $E$ denotes the classical expectation. Introducing the correct notion of extremals and variations we obtain
two different stochastic analogue of the least action principle depending on the regularity class we choose for
the admissible variations. The main point is that for variations in the class of Nelson differentiable process, the
extremals of our functional coincide with the stochastic Euler-Lagrange equation obtained via the stochastic embedding
procedure. This result is called the coherence lemma. In the reversible case, i.e. taking as a stochastic derivative only the
real part of our operator, we obtain the same result but in this case one can consider general variations.\\

In chapter \ref{noetherchap} we provide a first study of what dynamical data remain from the classical dynamical system under
the stochastic embedding procedure. We have focused on symmetries of the underlying equation and as a consequence on first
integrals. We prove a stochastic analogue of the Noether theorem. This allows us to define a natural notion of first
integral for stochastic differential equations. This part also put in evidence the need for a geometrical setting
governing Lagrangian systems which is the analogue of symplectic manifolds.\\

Chapter \ref{chapschrodinger} deals with the stochastic Euler-Lagrange equation for natural Lagrangian systems, i.e.
associated to Lagrangian functions of the form
\begin{equation}
L(x,v)=T(v)-U(x) ,
\end{equation}
where $U$ is a smooth function and $T$ is a quadratic form. In classical mechanics $U$ is the potential energy and $T$
the kinetic energy. The main result of this chapter is that by restricting our attention to good diffusion processes,
and up to a a well chosen function $\psi$, called the wave function, the stochastic Euler-Lagrange equation is equivalent
to a non linear Schr\"odinger equation. Moreover, by specializing the class of stochastic processes, we obtain the classical
Schr\"odinger equation. In that case, we can give a very interesting characterization of stochastic processes which are
solution of the stochastic Euler-Lagrange equation. Indeed, the square of the modulus of $\psi$ is equal to the density
of the associated stochastic process solution.\\

In chapter \ref{hamilton}, we define a natural notion of stochastic Hamiltonian system. This result can be seen as a first
attempt to put in evidence the stochastic analogue of a symplectic structure. We define a stochastic momentum process and
prove that, up to a suitable modification of the stochastic embedding procedure called the Hamiltonian stochastic embedding, and
reflecting the fact that the ``speed" of a given stochastic process is complex, we obtain a coherent picture with the classical
formalism of Hamiltonian systems. This first result is called the Legendre coherence lemma as it deals with the coherence
between the Hamiltonian stochastic embedding procedure and the Legendre transform. Secondly, we develop a Hamilton least
action principle and we prove again a coherence lemma, i.e. that the following diagram commutes
$$
\xymatrix{
  & H(x(t),p(t)) \ar[d]_{\mbox{\rm Hamilton least action principle}} \ar[r]^{S_H} & H(X(t),P(t)) \ar[d]^{
  \mbox{\rm Stochastic Hamilton least action principle}}       \\
  & (HE) \ar[r]_{S_H}   & (SHE)               }
$$
where $S_H$ denotes the Hamiltonian stochastic embedding procedure.\\

The last chapter discuss many possible developments of our theory from the point of view of mathematics and applications.

\part{The stochastic derivative}

\chapter{About Nelson stochastic calculus}
\label{chapnelson}

\section{About measurement and experiments}

In this section, we explain what we think are the basis of all possible
extensions of the classical derivative. The setting
of our discussion is the following:\\

We consider an experimental set-up which produces a dynamics. We assume that
each dynamics is observed during a time which is fixed, for example $[0,T]$,
where $T\in \R^{*+}$. For each experiment $i$, $i\in \N$, we denote by $X_i
(t)$ the dynamical
variable which is observed for $t\in [0,T]$.\\

Assume that we want to describe the {\it kinematic} of such a dynamical variable. What is the strategy ?\\

The usual idea is to {\it model} the dynamical behaviour of a variable by
ordinary differential equations or partial differential equations. In order to
do this, we must first try to have access to the {\it speed} of the variable. In order to
compute a significant quantity we can follow at least two different strategies:\\

\begin{itemize}
\item We do not have access to the variable $X_i (t)$, $t\in [0,T]$, but to a collection of measurements of this
dynamical variable. Assume that we want to compute the speed at time $t$. We can only compute an approximation of it for
a given resolution $h$ greater than a given threshold $h_0$. Assume that for each experiment we are able to compute the
quantity
\begin{equation}
v_{i,h} (t)=\di {X_i (t+h) -X_i (t) \over h} .
\end{equation}
We can then try to look for the behaviour of this quantity when $h$ varies. If
the underlying dynamics is not too irregular, then we can expect a limit for
$v_{i,h} (t)$ when $h$ goes to zero that we denote by $v_i (t)$.\\

We then compute the mean value
\begin{equation}
\bar{v} (t)=\di{1\over n} \di\sum_{i=1}^n v_i (t) .
\end{equation}
If the underlying dynamics is not too irregular then $\bar{v} (t)$ can be used to model the problem. In the contrary
the basic idea is to introduce a random variable.\\

Remark that due to the intrinsic limitation for $h$ we never have access to $v_i (t)$ so that this procedure can not be
implemented.

\item Another idea is to look directly for the quantity
\begin{equation}
\bar{v}_{h,n} (t)= \di{1\over n} \sum_{i=1}^n v_{i,h} (t) .
\end{equation}
Contrary to the previous case, if there exists a well defined mean value $\bar{v}_h (t)$ when $n$ goes to infinity then we can have a
as close as we want approximation. Indeed it suffices to do sufficiently many experiences. We then look for the limit of $\bar{v}_h (t)$ when $h$ goes to zero.
\end{itemize}

For regular dynamics these two procedures lead to the same result as all these quantities are well defined and converge
to the same quantity. This is not the case when we deal with highly irregular dynamics. In that case the second
procedure is easily implemented contrary to the first one. The only problem is that we loose the geometrical
meaning of the resulting limit quantity with respect to individual trajectories as one directly take a mean on all
trajectories before taking the limit in $h$.\\

This second alternative can be formalized using stochastic processes and leads to the Nelson backward and forward
derivatives that we define in the next section.\\

We have take the opportunity to discuss these notions because the previous remarks proves that one can not justify the form
of the Nelson derivatives using a geometrical argument like the non differentiability of trajectories for a Brownian motion.
This is however the argument used by E. Nelson (\cite{ne2},p.1080) in order to justify the fact that we need a substitute
for the classical derivative when studying Wiener processes. This misleadingly suggest that the forward and backward
derivative capture this non differentiability in their definition, which is not the case.

\section{The Nelson derivatives}

Let $X(t)$, $0\leq t\leq 1$ be $d$-dimensional continuous random process
defined on a probability space $(\Omega , {\cal A}, P)$, where ${\cal A}$ is the $\sigma$-algebra of all measurable
events and $P$ is a probability measure defined on $\cal A$. We denote by $I$ the open interval $(0,1)$.

\begin{defi}
The random process $X(t)$, $a\leq t\leq b$, is an SO-process if each $X(t)$
belongs to $L^1 (\Omega )$ and the mapping $t\rightarrow
X(t)$ from $\rR$ to $L^1 (\Omega )$ is continuous.
\end{defi}

Let ${\cal P}=\{ {\cal P}_t \}$ and ${\cal F}=\{ {\cal F}_t \}$
be an increasing and a decreasing family of sub-$\sigma$-algebras,
respectively, such that $X(t)$ is ${\cal F}_t$-measurable and ${\cal
P}_t$-measurable. In other words, ${\cal F}$ and ${\cal P}$ are two filtration
to which $X(t)$ is adapted. We let $E[\bullet \mid {\cal B} ]$ denote the
conditional expectation with respect to any sub-$\sigma$-algebra ${\cal B}
\subset {\cal A}$.

\begin{defi}
The random process $X(t)$, $a\leq t\leq b$, is an $S1$-process if it is an
SO-process such that
\begin{equation}
DX(t)=\lim_{h\rightarrow 0^+} E \left [ {X(t+h)-X(t)\over h } \mid {\cal P}_t \right ] ,
\end{equation}
and
\begin{equation}
D_* X(t)=\lim_{h\rightarrow 0^+} E \left [ {X(t)-X(t-h)\over h} \mid {\cal F}_t \right ],
\end{equation}
exist in $L^1 (\Omega )$ and the mappings $t\mapsto
DX(t)$ and $t\mapsto D_* X(t)$ are both continuous from $\rR$ to $L^1 (\Omega )$.
\end{defi}

\begin{defi}
The random process $X(t)$, $a\leq t\leq b$, is an S2-process if it is an
S1-process, and
\begin{equation}
\sigma^2 X(t)=\lim_{h\rightarrow 0^+} E \left [ {(X(t+h)-X(t))^2 \over h} \mid {\cal P}_t
\right ] ,
\end{equation}
and
\begin{equation}
\sigma_*^2 X(t)=\lim_{h\rightarrow 0^+} E\left [ {(X(t+h)-X(t))^2 \over h} \mid {\cal
F}_t \right ] ,
\end{equation}
exist in $L^1 (\Omega)$.
\end{defi}

\begin{defi}
We denote by $\Cu (I)$ the totality of S2-processes with continuous sample paths, such that $X(t)$, $DX(t)$ and
$D_* X(t)$, $a\leq t\leq b$, all lie in the Hilbert space $L^2 (\Omega)$ and are continuous
functions of $t$ in $L^2 (\Omega )$.

A completion of ${\cal C}^1 (I)$ in the norm
\begin{equation}
\parallel X\parallel=\sup_{t\in I} (\parallel X(t)\parallel_{L^2 (\Omega )} +\parallel DX(t)\parallel_{L^2 (\Omega )}
+\parallel D_* X(t)
\parallel_{L^2 (\Omega )} ),
\end{equation}
is also denoted by $\Cu (I)$, where $\parallel .\parallel_{L^2 (\Omega )}$ denotes the norm
of Hilbert space $L^2 (\Omega )$.
\end{defi}

\begin{rema}
The main point in the previous definitions for a forward and backward
derivative of a stochastic process, is that the forward and backward
filtration are fixed by the problem. As a consequence, we have not an
intrinsic quantity only related to the stochastic process. A possible
alternative definition is the following:

\begin{defi}
\label{alternative}
Let $X$ be a stochastic process, and $\sigma (X)$ (resp.
$\sigma_* (X)$) the forward (resp. backward) adapted filtration. We define
\begin{eqnarray}
\bab X(t)=\lim_{h\rightarrow 0^+} h^{-1} E[X(t+h)-X(t)\mid \sigma (X_s ,0\leq s\leq t) ],\\
\bab_* X (t)=\lim_{h\rightarrow 0^+} h^{-1} E[X(t)-X(t-h)\mid \sigma (X_s ,t\leq s\leq 1) ] .
\end{eqnarray}
\end{defi}

In this case, we obtain intrinsic quantities, only related to the stochastic
process. However, these new operators behave very badly from an algebraic view
point. Indeed, without stringent assumptions on stochastic processes, we do not
have linearity of $\bab$ or $\bab_*$.\\

This difficulty is not apparent as long as one restrict attention to a {\it single} stochastic process.
\end{rema}

\section{Good diffusion processes}

We introduce a special class of diffusion processes for which we can explicitly compute the
derivative $D$, $D_*$, $DD_*$, $D_*D$, $D^2$ and $D_*^2$.

\begin{defi}
We denote by $\Lambda_d$ the space of diffusion processes $X$ satisfying the following conditions:\\

\ni i- $X$ solves a stochastic differential equation :
\begin{eqnarray}
\label{SDE}
dX(t)=b(t,X(t))dt+\sigma(t,X(t))dW(t),\ \ X(0)=X_0 ,
\end{eqnarray}
where $X_0\in L^2(\Omega)$, $b:[0,T]\times \R^d\to\R^d$ and $\sigma:[0,T]\times
\R^d\to\R^d\otimes\R^d$ are Borel measurable functions satisfying the hypothesis
: there exists a constant $K$ such that for every $x,y\in\R^d$ we have
\begin{eqnarray}
    \sup_t \left(\left|\sigma(t,x)-\sigma(t,y)\right|+\left|b(t,x)-b(t,y)\right|\right)\leq K\left|x-y\right|,\\
    \sup_t \left(\left|\sigma(t,x)\right|+\left|b(t,x)\right|\right)\leq K(1+\left|x\right|).
\end{eqnarray}

\ni ii- For any $t>0$, $X(t)$ has a density $p_t(x)$ at point $x$.\\

\ni iii- Setting $a_{ij}=(\sigma\sigma^*)_{ij}$, for any $i\in\{1,\cdots,n\}$, for any $t_0>0$, for any bounded
open set $D\subset\R^d$,
\begin{equation}
    \int_{t_0}^1 \int_D \left|\partial_j(a_{ij}(t,x)p_t(x))\right|dxdt < +\infty.
\end{equation}

\ni iv- $b$ and $\d (t,x)\to \frac{1}{p_t(x)}\partial_j(a_{ij}(t,x)p_t(x))$ are continuous and bounded functions.
\end{defi}

\begin{rema}
\begin{itemize}
\item Hypothesis iii) ensures that (\ref{SDE}) has a unique $t-$conti\-nuous solution $X(t)$.

\item Hypothesis i), ii) and iii) allow to apply theorem $2.3$ p.$217$ in \cite{mns}.

\item We may wonder in which cases hypothesis ii) holds. Theorem $2.3.2$
$p.111$ of \cite{nua} gives the existence of a density for all $t>0$ under the
H\"ormander hypothesis which is involved by the stronger condition that the
matrix diffusion $\sigma\sigma^*$ is elliptic at any point $x$. A simple
example is given by a SDE where $b$ is a $\mathcal{C}^{\infty}(I\times \R^d)$
function with all its derivatives bounded, and where the diffusion matrix is a
constant equal to $c Id$. In this case, $p_t(x)$ belongs to
$\mathcal{C}^{\infty}(I\times \R^d)$; moreover, if $X_0$ has a differentiable
and everywhere positive density $p_0(x)$ with respect to Lebesgue measure such
that $p_0(x)$ and $p_0(x)^{-1}\nabla p_0(x)$ are bounded, then $b(t,x)-c \nabla
log(p_t(x))$ is bounded as noticed in the proof of proposition $4.1$ in
\cite{thieu}. So hypothesis ii) seems not to be such a restrictive
condition.

\item Assumption iv) is necessary to compute {\it explicitly} the second order operators of $D$ and $D_*$. The
existence of $D$ and $D_*$ is ensured under a weaker condition, the finite entropy condition equivalent to
\begin{equation}
E\left [
\int_0^1 (b(t,X(t))^2 \, dt \right ] <\infty .
\end{equation}
We refer to F\"ollmer (\cite{fol},proposition 2.5 p.121 and lemma 3.1 p.123) for more details.
\end{itemize}
\end{rema}

According to the theorem $2.3$ of \cite{mns} and thanks to iv), we
will see that $\Lambda_d\subset\EuScript{C}^1([0,T])$ and that we can compute $DX$ and
$D_*X$ for $X\in \Lambda_d$ (see Theorem \ref{derivgood}).

\section{The Nelson derivatives for good diffusion processes}

A useful property of good diffusions processes is that their Nelson's
derivatives can be explicitly computed. Precisely, we have:

\begin{thm}
\label{derivgood}
Let $X\in \Lambda_d$ which writes
$dX(t)=b(t,X(t))dt+\sigma(t,X(t))dW(t)$. Then $X$ is Markov diffusion with
respect to an increasing filtration $(\mathcal{P}_t)$ and a decreasing
filtration $(\mathcal{F}_t)$. Moreover, $DX$ and $D_*X$ exists w.r.t. these
filtration and :
\begin{eqnarray}
DX(t) & = & b(t,X(t))\\
D_*X(t) & = & b_*(t,X(t))
\end{eqnarray}
where $x\to p_t(x)$ denotes the density of $X(t)$ at $x$ and $$
b_*^i(t,x)=b^i(t,x)-\frac{1}{p_t(x)}\partial_j(a^{ij}(t,x)p_t(x))$$ with the
convention that the term involving $\frac{1}{p_t(x)}$ is $0$ if $p_t(x)=0$.
\end{thm}

\begin{proof}
The proof uses essentially theorem 2.3 of Millet-Nualart-Sanz \cite{mns} and the techniques of M. Thieullen for
the proof of proposition 4.1 in \cite{thieu}.\\

(1) Let $X\in \Lambda_d$. Then $X$ is a Markov diffusion w.r.t. the increasing
filtration $(\mathcal{P}_t)$ generated by the Brownian Motion $W(t)$ and so :
$$E\left[\frac{X(t+h)-X(t)}{h}\left|\mathcal{P}_t\right.\right]=E\left[\frac{1}{h}\int_t^{t+h}b(s,X(s))ds\left|\mathcal{P}_t\right.\right],
$$ and
$$E\left[\left|E\left[\frac{X(t+h)-X(t)}{h}\left|\mathcal{P}_t\right.\right]-b(t,X(t))\right|\right]\leq
E\left[\frac{1}{h}\int_t^{t+h}\left|b(s,X(s))-b(t,X(t))\right|ds\right].$$ We
can apply the dominated convergence theorem since $b$ is bounded and
$$\d \frac{1}{h}\int_t^{t+h}\left|b(s,X(s))-b(t,X(t))\right|ds \stackrel{h\to 0}{\longrightarrow} 0\ a.s.$$
(for $b$ is continuous and $X$ has a.s. continuous paths).\\

Therefore $DX$ exists and $DX(t)=b(t,X(t))$.\\

(2) As $X\in \Lambda_d$, we can apply theorem $2.3$ in \cite{mns}. So $\overline{X}(t)=X(1-t)$ is a diffusion process w.r.t. an increasing filtration $(\overline{\mathcal{P}}_t)$ and whose generator reads $\overline{L}_tf=\overline{b}^i\partial_i f+\frac{1}{2}\overline{a}^{ij}\partial_{ij}f$ with $\overline{a}^{ij}(1-t,x)=a^{ij}(t,x)$ and $\d \overline{b}^{i}(1-t,x)=-b^i(t,x)+\frac{1}{p_t(x)}\partial_j(a^{ij}(t,x)p_t(x))$.\\
Setting $\mathcal{F}_t=\overline{\mathcal{P}}_{1-t}$, $X$ is a Markov diffusion
w.r.t. the decreasing filtration $(\mathcal{F}_t)$. We have :
\begin{eqnarray}
    E\left[\frac{X(t)-X(t-h)}{h}\left|\mathcal{F}_t\right.\right] & = & E\left[\frac{\overline{X}(1-t)-\overline{X}(1-t+h)}{h}\left|\overline{\mathcal{P}}_{1-t}\right.\right] \nonumber\\
    & = & -E\left[\frac{1}{h}\int_{1-t}^{1-t+h}\overline{b}(s,\overline{X}(s))ds\left|\overline{\mathcal{P}}_{1-t}\right.\right].
\end{eqnarray}
Using the same calculations and arguments as above (since hypothesis iv) in the definition of class $\Lambda_d$ implies
that $\overline{b}$ is continuous and bounded), we obtain that $D_*X(t)$ exists and is equal to
$-\overline{b}(1-t,\overline{X}(1-t))$.
\end{proof}

In the case of {\it fractional Brownian motion} of order $H\not= 1/2$, the Nelson
derivatives do not exist. However, one can define new operators using the
so-called quasi conditional expectation introduced by \cite{aopu}. We refer to
the work of Darses and Sausserau \cite{ds} for more details.

\section{A remark about reversed processes}
\label{reversedresults}

This part reviews basic results about reversed processes, with a special
emphasis to diffusion processes. We use Nelson's
stochastic calculus.\\

Let $X$ be a process in the class $\mathcal{C}^1([0,1])$. We denote by
$\widetilde{X}$ the reversed process : $\widetilde{X}(t)=X(1-t)$, with his
"past" $\widetilde{\mathcal{P}}_t$ and his "future"
$\widetilde{\mathcal{F}}_t$.
As a consequence, we also have $\widetilde{x}\in\mathcal{C}^1([0,1]\to H)$.\\

Using the operators $\bab$ and $\bab_*$ defined in definition \ref{alternative},we
have:

\begin{lem}
$\bab_*x(t)=-\bab\widetilde{x}(1-t)=-\widetilde{\bab\widetilde{x}}(t).$
\end{lem}

\begin{proof}
The definition of $\bab_*$ gives immediately:
$$
\bab_*x(t)=\lim_{\epsilon\to
0^+}E\left[\left.\frac{\widetilde{x}(1-t)-\widetilde{x}(1-t+\epsilon)}{\epsilon}\right|\mathcal{F}_t\right].$$
But $\mathcal{F}_t=\sigma\{x(s),t\leq s\leq 1\}=\sigma\{\widetilde{x}(u),0\leq u\leq 1-t\}=\widetilde{\mathcal{P}}_{1-t}$.\\
Thus:
$$\bab_*x(t)=\lim_{\epsilon\to
0^+}-E\left[\frac{\widetilde{x}(1-t+\epsilon)-\widetilde{x}(1-t)}{\epsilon}\left|\widetilde{\mathcal{P}}_{1-t}\right.\right]
=-\bab\widetilde{x}(1-t) =-\widetilde{\bab\widetilde{x}}(t).$$
\end{proof}

The same computation is not at all possible when dealing with the operators $D$
and $D_*$.

\chapter{Stochastic derivative}
\label{chapstocder}
\setcounter{section}{0}
\setcounter{equation}{0}

In this part, we construct a natural {\it extension}\footnote{A precise meaning
to this word will be given in the following. It should be noted that Malliavin
calculus is not an extension of the ordinary differential calculus (see
below).} of the classical derivative on real stochastic processes as a unique
solution to an algebraic problem. This stochastic derivative turns out to be
necessarily complex valued. Our construction relies on Nelson's stochastic
calculus \cite{ne1}. We then study properties of our stochastic derivative and
establish a number of technical results, including a generalization of Nelson's
product rule \cite{ne1} as well as the stochastic derivative for functions of
diffusion processes . We also compute the stochastic derivative in some
classical examples. The main point is that, after a natural extension to
complex processes, the real part of the second derivative of a real stochastic
process coincide with Nelson's mean acceleration. We define a special class of
processes called Nelson differentiable, which will be of importance for the
stochastic calculus of variations developed in chapter \ref{chapvar}. This part is self
contained and all basic results about Nelson's stochastic calculus are
reminded.

\section{The abstract extension problem}

In this section, we discuss in a general abstract setting, what kind of
analogue of the classical derivative we are waiting for
on stochastic processes.\\

We first remark that real\footnote{Our aim was first to study dynamical systems
over $\R^n$. However, as we will see we will need to consider complex valued
objects.} valued functions naturally embed in stochastic processes.\\

Indeed, let $f:\R \rightarrow \R$ be a
given function. We denote by $X_f$ the deterministic stochastic process defined by
\begin{equation}
X_f (\omega )=f\ \forall \omega \in \Omega .
\end{equation}
We denote by $\iota : \R^{\R} \rightarrow {\cal P}$ the map associating to $f\in \R^{\R}$ the stochastic process $X_f$.\\

We denote by ${\cal P}_{\rm det}$ the subset of ${\cal P}$ consisting of deterministic processes, and by
${\cal P}^k_{\rm det}$ the set $\iota (C^k )$, $k\geq 1$.\\

As a consequence, we have a natural action of the classical derivative on the set of
differentiable deterministic processes, that we denote again $d/dt$.

Let $K=\R$ or $\C$. In the sequel, we denote by ${\cal P}_K \subset {\cal
S}_K$ a subset of the set of $K$-valued stochastic processes\footnote{We do not
give more precisions on this set for the moment, the set $\cal P$ can be the
whole set of real or complex valued stochastic processes, or a particular class
like diffusion processes,...etc.}.\\

Let $K=\R$ or $\C$.

\begin{defi}
Let $K=\R$ or $\C$. An extension of $d/dt$ on ${\cal P}_K$ is an operator
$\delta$, i.e. a map $\delta : {\cal P}_K \rightarrow {\cal S}_K$
such that:\\

{\rm i)} $\delta$ coincides with $d/dt$ on ${\cal P}^1_{\rm det}$,

{\rm ii)} $\delta$ is $\rR$-linear.
\end{defi}

Condition i), which is a gluing condition on the classical derivative is
necessary as long as one wants to relate classical
differential equations with their stochastic counterpart.\\

Condition ii) is more delicate. Of course, one has linearity of $\delta$ on
$\mbox{\rm Diff}$. A natural idea is then to preserve fundamental algebraic
properties of $d/dt$, $\R$-linearity being one of them. This condition is not
so stringent, if for example we consider $K=\C$. But, following this point of
view, one can ask for more precise properties like the Leibniz rule
\begin{equation}
\label{leib} d/dt (X\cdot Y)=d/dt (X) \cdot Y +X\cdot d/dt (Y) ,\ \ \forall X,Y\in {\cal P}^1_{\rm det} .
\end{equation}
In what follows, we construct a stochastic differential calculus based on
Nelson's derivatives.

\section{Stochastic differential calculus}

In this part, we extend the classical differential calculus to stochastic
processes using a previous work of Nelson \cite{ne1} on the dynamical theory of
Brownian motion. We define a stochastic derivative and review its properties.

\subsection{Reconstruction problem and extension}

Let us begin with some heuristic remarks supporting our definition and construction of a stochastic derivative.\\

Our aim is to construct a "natural" operator on $\Cu (I)$ which reduces to the
classical derivative $d/dt$ over differentiable deterministic
processes\footnote{A rigourous meaning to this sentence will be given in the
sequel.}. The basic idea underlying the whole construction is that, for
example in the case of the Brownian motion, the trajectories are
non-differentiable. At least, this is the reason why Nelson \cite{ne1}
introduces the left and right derivatives $DX$ and $D_* X$ for a given process
$X$. If we refer to geometry, forgetting for a moment processes for
trajectories, the fundamental property of the classical derivative $dx/dt (t_0
)$ of a trajectory $x(t)$ at point $t_0$, is to provide a first order
(geometric) approximation of the curve in a {\it neighbourhood} of $t_0$. One
wants to construct an operator, that we denote by $\cal D$, such that the data
of ${\cal D} X (t_0 )$ allows us to give an approximation of $X$ in a
neighbourhood of $t_0$. The difference is that we must know two quantities,
namely $DX$ and $D_* X$, in order to obtain the information\footnote{This
remark is only valid for general stochastic processes. Indeed, as we will see,
for diffusion processes, there is a close connection between $DX$ and $D_* X$,
which allows to simplify the definition of $\cal D$.}. For computational
reasons, one wants an operator with values in a field $F$. This field must be a
natural extension of $\rR$ (as we want to recover the classical derivative) and
at least of dimension $2$. The natural candidate to such a field is $\C$. One
can also recover $\C$ by saying that we must consider not only $\rR$ but the doubling algebra which corresponds to $\C$.\\

This informal discussion leads us to build a complex valued operator ${\cal
D} :\Cu (I) \rightarrow \Cu_{\cC} (I)$,
with the following constraints:\\

i) ({\it Gluing property}) For $X\in {\cal P}^1_{\rm det}$, ${\cal D}X(t)=dX/dt$,\\

ii) The operator $\cal D$ is $\R$-linear,\\

iii) ({\it Reconstruction property}) For $X\in \Cu (I)$, let us denote by
$${\cal D}X =A(DX ,D_* X)+iB(DX,D_* X),$$
where $A$ and $B$ are linear $\R$-valued mappings by ii). We assume that the
mapping
$$(DX ,D_* X) \mapsto (A(DX ,D_* X),B(DX ,D_* X))$$
is invertible.

\begin{lem}
\label{recons1}
The operator ${\cal D}$ has the form
$${\cal D}_{\mu} X= \left [ a DX +(1-a) D_* X\right ] +i \mu b\left [ DX -D_* X \right ]  ,\ \mu =\pm 1,$$
where $a,b\in \R$ and $b\not= 0$.
\end{lem}

\begin{proof}
We denote by $A(X)=a DX +b D_* X$ and $B(X)=c DX + d D_* X$. If $X\in C^1 (I)$,
we have $DX=D_* X=dX/dt$, and i) implies
$$a+b =1,\ c+d =0.$$
We then obtain the desired form. By iii), we must have $b\not= 0$ in order to
have invertibility.
\end{proof}

In order to rigidify this operator, we impose a constraint coming from the
analogy with the construction of the scale-derivative
for non-differentiable functions in \cite{cr1}.\\

iv) If $D_* =-D$, then $A(X)=0$, $B(X)=D$.\\

We then obtain the following result:

\begin{lem}
An operator $\cal D$ satisfying conditions i), ii), iii) and iv) is
of the form
\begin{equation}
{\cal D}_{\mu} =\di {D+D_* \over 2} +i\mu {D-D_* \over 2} ,\ \mu =\pm 1.
\end{equation}
\end{lem}

\begin{proof}
Using lemma \ref{recons1}, iii) implies the relations: $2a-1=0$ and $2b=1$, so
$a=b=1/2$.
\end{proof}

We then introduce the following notion of {\it stochastic derivative}:

\begin{defi}
We denote by ${\cal D}_{\mu}$ the operators defined by
$${\cal D}_{\mu} =\di {D+D_* \over 2} +i\mu {D-D_* \over 2} ,\ \mu =\pm 1.$$
\end{defi}

\subsection{Extension to complex processes}

In order to embed second order differential equations, we need to define the
meaning of ${\cal D}^2$, and more generally of ${\cal D}^n$, $n\in \nN$. The
basic problem is that, contrary to what happens for the ordinary differential
operator $d/dt$, even if we consider real valued processes $X$, the derivative
${\cal D}X$ is a complex one. As a consequence, one must extend
$\cal D$ to {\it complex processes}.\\

For the moment, let us denoted by ${\cal D}_{\cC}$ the extension to be define
of $\cal D$, to complex processes. Let $F$ be a field containing $\C$ to be
defined, and ${\cal D}_{\C} :\Cu_{\C} (I)\rightarrow F$. There are essentially
two possibilities to extend the stochastic derivative leading to the same
definition: an algebraic and an analytic one.

\subsubsection{Algebraic extension}

Let us assume that:\\

i)  the operator ${\cal D}_{\cC}$ is $\rR$-linear.\\

Let $Z=X +i Y$ be a complex process, where $X$ and $Y$ are two real processes.
By $\rR$-linearity, we have
$${\cal D}_{\cC} (Z)={\cal D}_{\cC} X +{\cal D}_{\cC}(iY) .$$
As ${\cal D}_{\cC}$ reduce to $\cal D$ on real processes, we obtain
$${\cal D}_{\cC} (Z)={\cal D} X +{\cal D}_{\cC}(iY) ,$$
which reduce the problem of the extension to find a suitable definition of
$\cal D$ on {\it purely imaginary
processes}.\\

We now make an assumption about the image of ${\cal D}_{\cC}$:\\

ii) The operator ${\cal D}_{\cC}$ is $\cC$-valued.\\

This assumption is far from being trivial, and has many consequences. One of
them is that, whatever the definition of ${\cal D}_{\cC} (iY)$ is, we will obtain a
complex quantity which {\it mixes} with the quantity ${\cal D} X$ in a non
trivial way.

\begin{rema}
One can wonder if another choice is possible, as for example, using {\it
quaternions} in order to avoid this {\it mixing problem}. However, a heuristic
idea behind the complex nature of $\cal D$ is that it corresponds to a
fundamental property of Nelson processes, the (in general) non-differentiable
character of trajectories. Then, the doubling of the underlying algebra is
related to a {\it symmetry breaking}\footnote{This reduces to $DX=D_* X$ for
deterministic differentiable processes, namely the invariance under
$h\rightarrow -h$.}. The computation of ${\cal D}^2$ is not related to such phenomenon.
\end{rema}

In the following, we give two different extensions of $\cal D$ to complex
processes under hypothesis i) and ii). The basic
problem is the following:\\

Let $Y$ be a real process. We denote
\begin{equation}
{\cal D} Y= S(Y)\pm i A(Y),
\end{equation}
where
\begin{equation}
S(Y)=\left [ \di {D+D_* \over 2} \right ] (Y), \ \mbox{\rm and}\ A(Y)=\left [ \di {D-D_* \over 2} \right ] (Y),
\end{equation}
and the letters $S$ and $A$ stand for the symmetric and antisymmetric operators
with respect to the exchange of $D$ with
$D_*$.\\

We denote
\begin{equation}
{\cal D}_{\cC} (iY)=R(Y)+iI(Y),
\end{equation}
where $R(Y)$ and $I(Y)$ are two real processes.\\

One can ask if we expect for special relations between $R(Y)$, $I(Y)$ and $S(Y)$,
$A(Y)$.

\paragraph{$\cC$-linearity}

If no relations are expected for, the natural hypothesis is to assume
$\cC$-linearity of ${\cal D}_{\cC}$, i.e.
\begin{equation}
{\cal D}_{\cC} (iY)=i{\cal D} Y.
\end{equation}
As a consequence, we obtain the following definition for the operator ${\cal D}_{\cC}$:\\

We denote by $\Cu_{\cC} (I)$ the set of stochastic processes of the form
$Z=X+iY$, with $X,Y\in \Cu (I)$.

\begin{defi}
The operator ${\cal D}_{\cC} :\Cu_{\cC} \rightarrow \Cu_{\cC}$ is defined by
$${\cal D}_{\cC ,\mu } (X+iY)={\cal D}_{\mu} X +i\mu {\cal D}_{\mu } Y ,\ \mu=\pm 1,$$
where $X,Y\in \Cu$.
\end{defi}

In the sequel, we denote ${\cal D}_{\cC}$ for ${\cal D}_{\cC ,\sigma}$.\\

The following lemma gives a strong reason to choose such a definition of ${\cal
D}_{\cC}$. We denote by
$${\cal D}_{\cC}^n ={\cal D}_{\cC} \circ \dots \circ {\cal D}_{\cC} .$$

\begin{lem}
We have
\begin{eqnarray}
\mathcal{D}_{\cC}^2=\left [{DD_*+D_*D \over 2} \right ] +i\left [{D^2-D_*^2
\over 2} \right ] .
\end{eqnarray}
\end{lem}

\begin{proof}
One use the $\C$-linearity of operator $\mathcal{D}$.
\end{proof}

We note that the real part of $\mathcal{D}^2$ is the {\it mean acceleration} as
defined by Nelson \cite{ne1}.

\begin{rema}
In (\cite{ne1},p.81-82), Nelson discusses natural candidates for the stochastic
analogue of acceleration. More or less, the idea
is to consider quadratic combinations of $D$ and $D^*$, respecting a gluing property with the classical derivative:\\

Let $Q_{a,b,c,d} (x,y)=ax^2 +bxy+cyx+dy^2$ be a real non-commutative quadratic
form such that $a+b+c+d=1$. A possible definition
for a stochastic acceleration is $Q(D,D^*)$.\\

We remark that the condition $a+b+c+d=1$ implies that when $D=D^*$, we have $Q(D,D^* )=D=D_*$.\\

The simplest examples of this kind are: $D^2$, $D_*^2$, $DD_*$ and $D_* D$.\\

We can also impose a {\it symmetry condition} in order to take into account
that we do not want to give a special importance to the mean-forward or
mean-backward derivative, by assuming that $Q(x,y)=Q(y,x)$, so that $Q$ is of
the form
$$Q_a (x,y)=a(x^2 +y^2) +(1-2a)\di {xy+yx\over 2}, a\in\R .$$
The simplest example in this case is obtained by taking $a=0$, i.e.
$$Q_0 (D,D_*)=\di {DD_* +D_* D \over 2} .$$
This last one corresponds to Nelson's mean acceleration and coincide with the real part of our stochastic derivative.\\

It must be pointed out that Nelson discuss only five possible candidates where
at least a three parameters family can be
defined by $Q_{a,b,c,1-a-b-c} (D,D_* )$. His five candidates correspond to the simplest cases we have described.\\

The choice of $Q_0 (D,D_*)$ as a mean acceleration is justified by Nelson using
a Gaussian Markov process $X(t)$ in equilibrium, satisfying the stochastic
differential equation
$$dX(t)=-\omega X(t)dt+dW(t).$$
We will return to this problem below.
\end{rema}

\subsubsection{Analytic extension}

We first remark that $D$ and $D_*$ possess a natural extension to complex
processes. Indeed, let $X=X_1 +iX_2$, with $X_i \in \Cu (I)$ then
$$D (X_1 +iX_2 )=D(X_1 )+iD(X_2 )\ \mbox{\rm and}\ D_* (X_1 +iX_2 )=D_* (X_1) +iD_* (X_2) .$$
As a consequence, the quantities $S(Y)$ and $A(Y)$ introduced in the previous
section for real valued processes make sense for complex processes, and the
quantity $A(X)+iS(X)$ is well defined for the complex process $X\in \Cu_{\C} (I)$.
As a consequence, we can naturally extend ${\cal D} (X)$ to complex processes
by simply posing
$${\cal D} (X)= \di {D+D_* \over 2} +\mu i{D-D_* \over 2}, $$
with the natural extension of $D$ and $D_*$.

\subsubsection{Symmetry}

A possible way to extend $\cal D$ is to assume that the regular part of ${\cal
D}_{\cC} (iY)$ is equal the imaginary part of ${\cal D} (Y)$, i.e. that the
geometric meaning of the complex and real part of ${\cal D}Y$ is exchanged. We
then impose the following relation:
$$R(Y)=\sigma A(Y).$$
This leads to the following extension:

\begin{defi}
The operator ${\cal D}_{\cC} :\Cu_{\cC} \rightarrow \Cu_{\cC}$ is defined by
$${\cal D}_{\cC ,\mu } (X+iY)={\cal D}_{\mu } X -i\mu {\cal D}_{\mu } Y ,\ \mu =\pm 1,$$
where $X,Y\in \Cu$.
\end{defi}

\subsection{Stochastic derivative for functions of diffusion process}

In the following, we need to compute the stochastic derivative of $f(t,X_t )$
where $X_t$ is a diffusion process and $f$ is a smooth function. Our main
result is the following lemma:

\begin{lem}
\label{darsien}
Let $X\in \Lambda_d$ and $f\in C^{1,2}(I\times \R^d)$ such that
$\partial_t f$, $\nabla f$ and $\partial_{ij}f$ are bounded. Then, we have:
\begin{eqnarray}
    Df(t,X(t)) & = & \left[\partial_t f + DX(t)\cdot \nabla f +\frac{1}{2}a^{ij}\partial_{ij}f\right](t,X(t)) ,\\
    D_*f(t,X(t)) & = & \left[\partial_t f + D_*X(t)\cdot \nabla f -\frac{1}{2}a^{ij}\partial_{ij}f\right](t,X(t)) .
\end{eqnarray}
\end{lem}

\begin{proof}
Let $X\in \Lambda_d$ and $f\in C^{1,2}(I\times \R^d)$ such that $\partial_t f$, $\nabla
f$ and $\partial_{ij}f$ are bounded. Thus $f$ belongs to the domain of the
generators $L_t$ and $\overline{L}_t$ of the diffusions $X(t)$ and
$\overline{X}(t)$. Moreover these regularity assumptions allow us to use the
same arguments as in the proof of theorem \ref{derivgood} in order to write :
\begin{eqnarray}
Df(t,X(t)) & = & \partial_t f(t,X(t)) +L_t(f(t,\cdot))(X(t))\nonumber\\
                     & = & \left[\partial_t f+b^i\partial_i f+\frac{1}{2}a^{ij}\partial_{ij}f\right](t,X(t))\nonumber\\
                     & = & \left[\partial_t f+DX(t)\cdot\nabla f+\frac{1}{2}a^{ij}\partial_{ij}f\right](t,X(t))\nonumber
\end{eqnarray}
and
\begin{eqnarray}
D_*f(t,X(t)) & = & \partial_t f(t,X(t)) -\overline{L}_{1-t}(f(t,\cdot))(X(t))\nonumber\\
                     & = & \left[\partial_t f+D_*X(t)\cdot\nabla f-\frac{1}{2}a^{ij}\partial_{ij}f\right](t,X(t))\nonumber
\end{eqnarray}
\end{proof}

We deduce immediately the following corollary :

\begin{cor}
Let $X\in \Lambda_d$ and $f\in C^{1,2}(I\times \R^d)$ such that $\partial_t f$, $\nabla
f$ and $\partial_{ij}f$ are bounded. Then, we have:
\begin{eqnarray}
\mathcal{D}_\mu f(t,X(t)) & = & \left[\partial_t f + \mathcal{D}_\mu X(t)\cdot
\nabla f +\frac{i\mu}{2}a^{ij}\partial_{ij}f\right](t,X(t)) .
\end{eqnarray}
\end{cor}

and

\begin{cor}
\label{corcomm}
Let $X\in \Lambda_d$ with a constant diffusion coefficient $\sigma$ and
$f\in C^{1,2}(I\times \R^d)$ such that $\partial_t f$, $\nabla f$ and
$\partial_{ij}f$ are bounded. Then, we have:
\begin{eqnarray}
\mathcal{D}_\mu f(t,X(t)) & = & \left[\partial_t f + \mathcal{D}_\mu X(t)\cdot
\nabla f +\frac{i\mu\sigma^2}{2}\Delta f\right](t,X(t)) .
\end{eqnarray}
\end{cor}

\subsection{Examples}

We compute the stochastic derivative in some famous examples, like the
Ornstein-Uhlenbeck process and a Brownian mation in an external force.

\subsubsection{The Ornstein-Uhlenbeck process}

A good model of the Brownian motion of a particle with friction is provided by
the Ornstein-Uhlenbeck equation:
\begin{equation}
    \left\{
\begin{array}{l}
    X''(t)=-\alpha X'(t)+\sigma\xi(t) \\
    X(0)=X_0,\ X'(0)=V_0,
\end{array}\right.
\end{equation}
where X(t) is the position of the particle at time, $\alpha$ is the friction
coefficient, $\sigma$ is the diffusion coefficient, $X_0$ and $V_0$ are given
Gaussian variables, $\xi$ is "white noise". The term $-\alpha X'(t)$ represents
a
frictional damping term.\\

The stochastic differential equation satisfied by the velocity process
$V(t):=Y'(t)$ is given by:
\begin{equation}
    \left\{
\begin{array}{l}
    dV(t)=-\alpha V(t)dt+\sigma dW(t) \\
    V(0)=V_0,
\end{array}\right.
\end{equation}

We can explicitly compute $\mathcal{D}V$ and $\mathcal{D}^2 V$:

\begin{lem}
Let $V(\cdot)$ be a solution of \begin{equation}
    \left\{
\begin{array}{l}
    dV(t)=-\alpha V(t)dt+\sigma dW(t) \\
    V(0)=V_0,
\end{array}\right.
\end{equation} where $V_0$ has a normal distribution with mean zero and variance
$\frac{\sigma^2}{2\alpha}$.\\

Then $V\in \EuScript{C}^2(]0,+\infty))$ and:
\begin{eqnarray}
    \mathcal{D}V(t) & = & -i\alpha V(t)\\
    \mathcal{D}^2V(t) & = & -\alpha^2 V(t).
\end{eqnarray}
\end{lem}

\begin{proof}
The solution is a Gaussian process explicitly given by:
\begin{equation}
    \forall t\geq 0,\ V(t)=V_0 e^{-\alpha t} + \sigma \int_0^t e^{-\alpha (t-s)}dW(s).
\end{equation}

Therefore, we can compute the expectation and the variance of the normal
variable $V(t)$ :
\begin{equation}
    \left\{
\begin{array}{l}
    E[V(t)]=E[V_0]e^{-\alpha t} \\
    {\rm Var}(V(t))=\frac{\sigma^2}{2\alpha}+\left({\rm Var} (V_0))-\frac{\sigma^2}{2\alpha}\right)e^{-2\alpha t},
\end{array}\right.
\end{equation}
We notice, as in \cite{karatzas}, that if $V_0$ has a normal distribution with
mean zero and variance $\frac{\sigma^2}{2\alpha}$, then $X$ is a stationary
gaussian process which distribution $p_t(x)$ at each time $t$ reads
\begin{equation}
p_t(x)=\frac{\sqrt{\alpha}}{\sqrt{\pi}\sigma}e^{-\frac{\alpha x^2}{\sigma^2}} .
\end{equation}
As a consequence, we have
\begin{equation}
\d \forall t\geq 0,
\ln(p_t(x))=\ln(\frac{\sqrt{\alpha}}{\sqrt{\pi}\sigma})-\frac{\alpha
x^2}{\sigma^2} ,
\end{equation}
and
\begin{equation}
\d \sigma^2 \partial_x \ln(p_t(x))=\sigma^2\frac{-2\alpha x}{\sigma^2}=-2\alpha
x .
\end{equation}
Moreover, we have
\begin{equation}
DV(t)=-\alpha V(t),
\end{equation}
and according to theorem \ref{derivgood}, we obtain
\begin{equation}
D_*V(t)=-\alpha V(t)-\sigma^2 \partial_x \ln(p_t(V(t)))=\alpha V(t).
\end{equation}
Therefore $\mathcal{D}V(t)=-i\alpha V(t)$, and using the $\C-$linearity of
$\mathcal{D}$, we obtain $\mathcal{D}^2V(t)=-\alpha^2 V(t)$, which concludes
the proof.
\end{proof}

\subsubsection{Brownian particle submitted to an external force}

In some examples of random mechanics, one has to consider the stochastic
differential system:
\begin{equation}
    \left\{
\begin{array}{l}
    dX(t)=V(t)dt\\
    dV(t)=-\alpha V(t)dt+K(X(t))dt+\sigma dW(t) \\
    X(0)=X_0,\ V(0)=V_0,
\end{array}\right.
\end{equation}
$X$ and $V$ may represent the position and the velocity of a particle of mass
$m$ being under the influence of an external force $F=-\nabla U$ where $U$ is a
potential. Set $K=F/m$.
The "free" case $K=0$ is the above example.\\

When $K(x)=-\omega^2 x$ (a linear restoring force), the system can also be seen
as the random harmonic oscillator. In this case, it can be shown that if
$(X_0,V_0)$ has an appropriate gaussian distribution then $(X(t),V(t))$ is
a stationary gaussian process in the same way as before.\\

Let us come back to the general case.\\

First, we remark that $X$ is Nelson-differentiable and we have
$DX(t)=D_*X(t)=V(t)$. Moreover, Nelson claims in (\cite{ne1},p.83-84) that,
when the particle is in equilibrium with a special stationary density,
\begin{eqnarray}
DV(t) & = & -\alpha V(t)+K(X(t)),\\
D_*V(t) & = & \alpha V(t) +K(X(t)).
\end{eqnarray}
We can summarize these results with the computation of $\mathcal{D}$ :

\begin{eqnarray}
    \mathcal{D}X(t) & = & V(t),\\
    \mathcal{D}^2 X(t) & = & K(X(t))-i\alpha V(t).
\end{eqnarray}

\newpage
\chapter{Properties of the stochastic derivatives}
\label{chapprop}
\setcounter{section}{0}
\setcounter{equation}{0}

\section{Product rules}

In chapter \ref{chapvar}, we develop a stochastic calculus of variations. In many problems,we
will need the analogue of the classical formula of {\it integration by parts},
based on the following identity, called the {\it product or Leibniz rule}
$$\di {d\over dt} (fg)=\di {df\over dt} g+f\di {dg\over dt}, \eqno{(P)}$$
where $f,g$ are two given functions.\\

Using a previous work of Nelson \cite{ne1}, we generalize this formula for our
stochastic derivative. We begin by recalling the fundamental result of Nelson
on a product rule formula for backward and forward derivatives:

\begin{thm}
Let $X,Y\in \Cu (I)$, then we have:
\begin{eqnarray}
\frac{d}{dt}E[X(t)\cdot Y(t)]=E[DX(t)\cdot Y(t)+X(t)\cdot D_*Y(t)]
\end{eqnarray}
\end{thm}

We refer to (\cite{ne1},p.80-81) for a proof.

\begin{rema}
It must be pointed out that this formula mixes the backward and forward
derivatives. As a consequence, even without our definition of the stochastic
derivative, which takes into account these two quantities, the previous product
rule suggests the construction of an operator which mixes these two terms in a
"symmetrical" way.
\end{rema}

We now take up the various consequences of this formula regarding our operator
$\mathcal{D}$. A straightforward calculation gives:

\begin{lem}
\label{lem_rule}
Let $X,Y\in \Cu(I)$, we then have:
\begin{eqnarray}
\frac{d}{dt}E[X(t)\cdot Y(t)] & = & E[\rea (\mathcal{D}X(t))\cdot Y(t)+X(t)\cdot \rea (\mathcal{D}Y(t))]    \\
E[\ima (\mathcal{D}X(t))\cdot Y(t)] & = & E[X(t)\cdot \ima (\mathcal{D}Y(t))]
\end{eqnarray}
\end{lem}

\begin{lem}
Let $X,Y\in\Cu_{\C}(I)$. We write $X=X_1+iX_2$ and $Y=Y_1+iY_2$ where
$X_i,Y_i\in\Cu (I)$. Therefore :
\begin{equation}
\label{formule} E[\mathcal{D}_{\mu}X\cdot Y+X\cdot \mathcal{D}_{\mu}Y] =
\frac{d}{dt}g(X(t),Y(t))+r(X(t),Y(t)) ,
\end{equation}
where
\begin{equation}
g(X,Y)=E[X\cdot Y] ,
\end{equation}
and
\begin{equation}
\left .
\begin{array}{lll}
r(X,Y) & = & -2E[Y_1\cdot \ima(\mathcal{D}_{\mu}X_2)]-2E[Y_2\cdot \ima(\mathcal{D}_{\mu}X_1)] \\
 & & +i\left( 2E[Y_1\cdot \ima(\mathcal{D}_{\mu}X_1)]  -2E[Y_2\cdot \ima(\mathcal{D}_{\mu}X_2)] \right) .
\end{array}
\right .
\end{equation}
\end{lem}

\begin{proof}
We have
\begin{equation}
\left .
\begin{array}{lll}
\d Y\mathcal{D}_{\mu} X & = &
Y_1\rea(\mathcal{D}_{\mu}X_1)-Y_1\ima (\mathcal{D}_{\mu}X_2) \\
 &  & -Y_2\ima (\mathcal{D}_{\mu}X_1) -Y_2\rea (\mathcal{D}_{\mu}X_2)\\
 &  & +i \left( Y_1\ima(\mathcal{D}_{\mu}X_1)
+Y_1\rea (\mathcal{D}_{\mu}X_2) +Y_2\rea (\mathcal{D}_{\mu}X_1)-Y_2\ima
(\mathcal{D}_{\mu}X_2) \right) .
\end{array}
\right .
\end{equation}
In a symmetrical way, we obtain
\begin{equation}
\left .
\begin{array}{lll}
\d X\mathcal{D}_{\mu}Y & = &
X_1\rea(\mathcal{D}_{\mu}Y_1)-X_1\ima (\mathcal{D}_{\mu}Y_2)\\
 & & -X_2 \ima(\mathcal{D}_{\mu}Y_1)-X_2\rea (\mathcal{D}_{\mu}Y_2)\\
 & & +i\left(X_1\ima (\mathcal{D}_{\mu}Y_1)
+X_1\rea (\mathcal{D}_{\mu}Y_2)+X_2\rea (\mathcal{D}_{\mu}Y_1)-X_2\ima
(\mathcal{D}_{\mu}Y_2)\right).
\end{array}
\right .
\end{equation}
Forming the sum of these expressions and using lemma \ref{lem_rule}, we obtain (\ref{formule}).
\end{proof}

The next lemma will be of importance in chapter \ref{chapvar} for the derivation of the
stochastic analogue of the Euler-Lagrange equations:

\begin{lem}
Let $X,Y\in\EuScript{C}^1_{\C}(I)$. We write $X=X_1+iX_2$ and $Y=Y_1+iY_2$
where $X_i,Y_i\in\EuScript{C}^1(I)$. Therefore, we have:
\begin{equation}
\label{formule2}
E[\mathcal{D}_{\mu}X\cdot Y+X\cdot \mathcal{D}_{-\mu}Y] = \frac{d}{dt}g(X(t),Y(t))
\end{equation}
where $g(X,Y)=E[X_1\cdot Y_1-X_2\cdot Y_2]+i E[Y_1\cdot X_2+Y_2\cdot
X_1]=E[X\cdot Y]$
\end{lem}

\begin{proof}
We have
\begin{equation}
\left .
\begin{array}{lll}
\d Y\mathcal{D}_{\mu}X & = & Y_1\Re(\mathcal{D}_{\mu}X_1) -Y_1\Im(\mathcal{D}_{\mu}X_2)\\
 & & -Y_2\Im(\mathcal{D}_{\mu}X_1)
-Y_2\Re(\mathcal{D}_{\mu}X_2)\\
 & & +i\left(Y_1\Im(\mathcal{D}_{\mu}X_1)
+Y_1\Re(\mathcal{D}_{\mu}X_2)+Y_2\Re(\mathcal{D}_{\mu}X_1)-Y_2\Im(\mathcal{D}_{\mu}X_2)\right) ,
\end{array}
\right .
\end{equation}
and in a symmetrical way
\begin{equation}
\left .
\begin{array}{lll}
\d X\mathcal{D}_{-\mu}Y & = & (X_1+iX_2)(\overline{\mathcal{D}_{\mu}Y_1}
+i\overline{\mathcal{D}_{\mu}Y_2}) \\
 & = & X_1\Re(\mathcal{D}_{\mu}Y_1)+X_1\Im(\mathcal{D}_{\mu}Y_2) \\
 &   & +X_2\Im(\mathcal{D}_{\mu}Y_1)-X_2\Re(\mathcal{D}_{\mu}Y_2) \\
 &   & +i\left(-X_1\Im(\mathcal{D}_{\mu}Y_1)+X_1\Re(\mathcal{D}_{\mu}Y_2)
 +X_2\Re(\mathcal{D}_{\mu}Y_1) +X_2\Im(\mathcal{D}_{\mu}Y_2)\right).
\end{array}
\right .
\end{equation}

We form the sum of these expressions and we use the lemma \ref{lem_rule} to
obtain (\ref{formule}).
\end{proof}

\subsection{A new algebraic structure}

A convenient way to write equation (\ref{formule2}) is to use the following Hermitian product:\\

For all $X,Y\in {\cal P}_{\cC}$, we denote by $\star$ the product
\begin{equation}
X\star Y=X\cdot \overline{Y} ,
\end{equation}
where $.$ denotes the usual scalar product.\\

Formula (\ref{formule2}) is then equivalent to:
\begin{equation}
{\cal D} E[X\star Y]= E\left [ {\cal D} X\star Y +X\star {\cal D} Y \right ] ,
\end{equation}
where we have implicitly used the fact that ${\cal D}$ reduces to $d/dt$ when this quantity has a sense.\\

This new form leads us to the introduction of the following algebraic structure, which
is, as far as we know, new. Let $\delta$ be the canonical mapping
\begin{equation}
\delta : \left .
\begin{array}{lll}
{\cal P}_{\cC} \otimes {\cal P}_{\cC} & \rightarrow & {\cal P}_{\cC}\\
X \otimes Y & \mapsto & X\star \overline{Y} .
\end{array}
\right .
\end{equation}
We define for ${\cal D}$ the quantity $\Delta ({\cal D} )={\cal D}\otimes 1 +1
\otimes {\cal D}$, which we will call the {\it coproduct} of ${\cal D}$. Then,
denoting by $E$ the classical mapping which takes the expectation of a given
stochastic process, we obtain the following diagram:

\begin{equation}
\begin{CD}
{\cal P}_{\cC} \otimes {\cal P}_{\cC}  @> \Delta (D) >> {\cal P}_{\cC} \otimes {\cal P}_{\cC}\\
X\otimes Y @> >> {\cal D} X \otimes Y + X\otimes {\cal D} Y \\
@VV{\delta}V @VV{\delta}V\\
X\star Y @>  >> {\cal D} X \star Y +X\star {\cal D} Y \\
@VV{\mbox{\rm E}}V @VV{\mbox{\rm E}}V\\
E[X\star Y] @> {\cal D} >> E [{\cal D} X \star Y +X\star {\cal D} Y]
\end{CD}
\end{equation}
This structure is similar to the classical algebraic structure of {\it Hopf
algebra}. The difference is that we perturb the classical relations by a linear
mapping, here given by $E$. It will be interesting to study this kind of
structure in full generality.

\section{Nelson differentiable processes}

\subsection{Definition}

We define a special class of processes, called Nelson-differentia\-ble processes,
which will play an important role in the stochastic calculus of variations of
chapter \ref{chapvar}.

\begin{defi}
A process $X\in \Cu (I)$ is called Nelson differentiable if $DX=D_* X$.
\end{defi}

\begin{notation}
We denote by ${\cal N}^1 (I)$ the set of Nelson differentiable processes.
\end{notation}

A better definition is perhaps to use $\cal D$ instead of $D$ and $D_*$ saying
that Nelson differentiable processes have a
real stochastic derivative.\\

The main idea behind this definition is that we want to define a class ${\cal
P}$ of processes in $\Cu (I)$ such that if $X\in \Cu (I)$ then for all $Y\in
{\cal P}$, we have
$$\ima ({\cal D}(X+Y))=\ima ({\cal D}X) .$$
This condition imposes that $\ima ({\cal D} Y)=0$.\\

This condition will appear more clearly in chapter \ref{chapvar} concerning the stochastic calculus of variations.

\begin{rema}
We must keep in mind that our definition of the stochastic derivative follows
the idea of the {\it scale calculus} developed in \cite{cr1} to study non-differentiable functions.
In that context, the existence of an imaginary part
for the scale derivative of a function is seen as a resurgence of its non-differentiability. In particular, when
the underlying function is differentiable
then the scale derivative is real. That is why we have chosen to call processes
such that $D=D_*$ Nelson differentiable.
\end{rema}

The definition of Nelson differentiable processes is only given for processes
in $\Cu (I)$. It is not at all clear to know what is the correct extension to
$\Cu_{\C} (I)$. As we have no use of such kind of notion on $\Cu_{\C} (I)$ we
don't discuss this point here.\\

Of course a difficult problem is to characterize these processes. The next section discusses some examples.

\subsection{Examples of Nelson-differentiable process}

We give examples of Nelson-differentiable processes.

\subsubsection{Differentiable deterministic process}

It is probably the first and the simplest example. Let $x(\cdot)$ be a
differentiable deterministic process defined on $I\times \Omega$. The past
$\mathcal{P}$ and the future $\mathcal{F}$ are trivial:
$$\forall t\in I,\ \mathcal{P}_t=\mathcal{F}_t=\{\emptyset,\Omega\}.$$
As a consequence, we have
$$\forall t\in I,\ Dx(t)=D_*x(t)=x'(t),$$
where $x'$ is the usual derivative of $x$.

\subsubsection{A very special random example}

Let $X\in\EuScript{C}^1(I)$. In \cite{ne1}, Nelson shows that $X$ is a constant
(i.e. X(t) is the same random variable for all t) if and only if : $\forall
t\in I,\ DX(t)=D_*X(t)=0$. So it provides us a random example of
$\mathcal{N}^1(I)-$process.

\subsubsection{Nelson-differentiable diffusion processes}

Using theorem \ref{derivgood}, we can find a sufficient and necessary condition
for a diffusion process to be a Nelson-differentiable process:

\begin{lem}
Let $X\in \Lambda_d$ with $\sigma=const$, then $X\in {\cal N}^1 (I)$ if and only if
\begin{equation}
\nabla(\sigma^2 p)(t,X(t))=0 .
\end{equation}
\end{lem}

When the diffusion equation is time homogeneous and the solutions have a
density, we note that this density must be a stationary density. Moreover, the
Fokker-Planck equation (Kolmogorov forward equation) allows us to give a necessary
condition (a relation between the drift and the diffusion coefficient) for a
diffusion equation to give a Nelson-differentiable solution.

\subsubsection{The random harmonic oscillator}

The random harmonic oscillator satisfies the stochastic differential equation:
\begin{equation}
    \left\{
\begin{array}{l}
    dX(t)=V(t)dt\\
    dV(t)=-\alpha V(t)dt-\omega^2 X(t)dt+\sigma dW(t) \\
    X(0)=X_0,\ V(0)=V_0,
\end{array}\right.
\end{equation}

As a consequence, we have $\d X(t)=\int_0^t V(s)ds$ with $\d E\left[\int_0^b
\left|V(s)\right|^2ds\right]<\infty$ ($b>0$), and $X$ has a strong derivative
in $L^2$. We then obtain $DX(t)=D_*X(t)=V(t)$. Finally, we have $X\in
\mathcal{N}^1([0,b])$ and $\mathcal{D}X(t)=V(t)$.

\subsection{Product rule and Nelson-differentiable processes}

\begin{cor}
Let $X,Y\in\EuScript{C}^1_{\C}(I)$. If $X$ is Nelson-differentiable then :
\begin{equation}
\label{derivation} E[\mathcal{D}_{\mu}X(t)\cdot Y(t)+X(t)\cdot
\mathcal{D}_{\mu}Y(t)]=\frac{d}{dt} E(X(t),Y(t))
\end{equation}
\end{cor}

\begin{proof}
This is a simple consequence of the fact that if $X=X_1+iX_2$ is
Nelson-differentiable then
$\mbox{\rm Im}(\mathcal{D}_{\mu}X_1)=\mbox{\rm Im}(\mathcal{D}_{\mu}X_2)=0$.
\end{proof}

\part{Stochastic embedding procedures}

\chapter{Stochastic embedding of differential operators}
\label{chapembed}
\setcounter{section}{0}
\setcounter{equation}{0}

A natural question concerning ordinary and partial differential equations concerns their behaviour
under small random perturbations.
This problem is particularly important in natural phenomena where we know that
models are only an approximation of the real setting. For example, the study of
the long term behaviour of the solar system is usually done by running numerical
computations on the $n$-body problem. However, many effects in the solar
systems are not included in this model and can be of importance if one looks for
a long term integration, as non conservative effects (due to tidal forces
between planets) and the oblatness of the sun which is not yet modelled by a
differential equation.

The main problem is then to find the correct analogue of a given differential
equation taking into account the following facts:\\

i) The classical equation is a good model at least in first approximation,

ii) One must extend this equation to stochastic processes.\\

Using the stochastic derivative introduced in the previous part, we give a
natural embedding of partial or ordinary differential equations into stochastic
partial or ordinary differential equations. It must be pointed out that we
do not perturb the classical equation by a random noise or anything else. In
this respect we are far from the usual way of thinking underlying the fields of
stochastic differential equations or stochastic dynamical systems.

Of course, having this natural embedding, we can naturally define what a stochastic perturbation of a differential
equation is. This is simply a stochastic
perturbation of the stochastic embedding of the given equation. The main point
is that we stay in the same class of objects dealing with perturbations, which
is not the case in the stochastic theory of differential equations, where we
jump from classical solutions to stochastic processes in one step using for
example Ito's stochastic calculus\footnote{This remark is also valid for all the
theories of this kind, using your
favourite stochastic calculus, like Malliavin calculus for example.}.\\

In this part we first give a general embedding procedure for partial
differential equations. We discuss classical examples, in particular first and
second order differential equations. The case of Lagrangian systems is studied
in details in chapter \ref{chapvar}. An important part of classical differential equations
coming from mechanics are {\it reversible}. This property is not conserved by
the previous stochastic embedding procedure. We define a special embedding
called reversible, which preserves this property, meaning that if $X$ is a
solution of the stochastic embedded equation, then $\tilde{X}$, the reversed
process, is again a solution.

\section{Stochastic embedding of differential operators}

In this part, we first give an abstract embedding procedure based on an
extension of the classical derivative defined in the previous part. We then
specialize our embedding procedure using the stochastic derivative.

\subsection{Abstract embedding}

Let $\kk$ be a ring, we denote by $\kk [x]$ the ring of polynomials with
coefficients in $\kk$. Let $\kk =C^1 (\rR^d \times \R)$.

\begin{defi}
A {\it differential operator} is an elements of $\kk [d/dt ]$.
\end{defi}

Let $O\in \kk [d/dt ]$, the differential operator $O$ is of the form
\begin{equation}
O=a_0 (\bullet,t)+a_1 (\bullet ,t)\di {d\over dt} +\dots + a_n (\bullet ,t) \di {d^n \over dt^n },
\ a_i \in \kk,\ =0,\dots ,n,
\end{equation}
for a given $n\in \nN$, called the {\it degree} of $O$.\\

The action of $O$ on a given function $x:\R \rightarrow \R^d$, $t\mapsto x(t)$ is denoted $O\cdot x$ and defined by
\begin{equation}
O\cdot x =\di\sum_{i=0}^n  a_i (x(t),t) \di {dx\over dt} .
\end{equation}

\begin{defi}[Abstract stochastization]
\label{stoca}
Let $O\in \kk [d/dt ]$ be a differential operator, of the form
\begin{equation}
O=a_0 (\bullet ,t)+a_1 (\bullet ,t)\di {d\over dt} +\dots + a_n (\bullet ,t)\di {d^n \over dt^n },
\ a_i \in \kk,\ =0,\dots ,n,
\end{equation}
where $n\in \nN$ is given.

The stochastic embedding of $O$ with respect to the extension $\delta :{\cal
P}\rightarrow {\cal P}$ is an element $O_{\delta}$ of ${\cal P}[\delta ]$
defined by
\begin{equation}
O_{\delta} =a_0 (\bullet ,t) +a_1 (\bullet ,t)\di \delta +\dots + a_n (\bullet ,t) \di {\delta^n }, \ a_i
\in {\cal P},\ i=0,\dots ,n,
\end{equation}
where $\delta^n =\delta \circ \dots \circ \delta$.\\

The action of $O_{\delta}$ on a given stochastic process $X$, denoted by $O_{\delta} \cdot X$ is defined by
\begin{equation}
O_{\delta} \cdot X =\di\sum_{i=0}^n a_i (X,t) \di\delta ^i X ,
\end{equation}
where the notation $a_i (X,t)$ stands for the stochastic process defined for all $\omega\in \Omega$ by
\begin{equation}
a_i (X,y) (\omega )=a_i (X(\omega,t), t) .
\end{equation}
\end{defi}

The main property of this embedding is the fact that
\begin{equation}
O_{\delta} \mid_{{\cal P}_{\rm det}^n} =O ,
\end{equation}
so that the classical differential equation associated to $O$, and
given by
$$
O\cdot x=0,\eqno{(E)}
$$
is contained in the stochastic differential equation
$$
O_{\delta} \cdot X=0 .\eqno{(SE)}.
$$

\subsection{Nelson Stochastic embedding}

Using the stochastic derivative, we have a particular stochastic embedding
procedure.

\begin{defi}[Stochastization]
\label{stoca}
Let $O\in \kk [d/dt ]$ be a differential operator, of the form
\begin{equation}
\label{formop}
O=a_0 (\bullet ,t)+a_1 (\bullet,t)\di {d\over dt} +\dots + a_n (\bullet,t)\di {d^n \over dt^n },
\ a_i \in \kk,\ =0,\dots ,n,
\end{equation}
where $n\in \nN$ is given.

The stochastic embedding of $O$ with respect to the stochastic extension ${\cal D}_{\mu}$ is an element $O_{\rm stoc}$ of
$C^1 (I)[{\cal D}_{\sigma} ]$
defined by
\begin{equation}
O_{\rm stoc} =a_0 (\bullet ,t) +a_1 (\bullet ,t)\di {\cal D} +\dots + a_n (\bullet ,t) \di {{\cal D}^n }, \ a_i
\in C^1 (I),\ i=0,\dots ,n.
\end{equation}
\end{defi}

We denote by $\stoc$ the operator associating to an operator $O$ of the form \ref{formop} the operator $O_{\rm stoc}$.
As a consequence, we will frequently use the notation $\stoc (O)$ for $O_{\rm stoc}$.\\

In some occasions, in particular for the Euler-Lagrange equation, we will need to consider differential operators in a
non-standard form. Precisely, we need to consider operators like
\begin{equation}
\di B_a = {d\over dt} \circ a (\bullet ,t) .
\end{equation}
This notation means that $B_a$ acts on a given function as
\begin{equation}
B_a \cdot x =\di {d\over dt} \left ( a(x(t),t)) \right ) .
\end{equation}
The basic idea is to define the stochastic embedding of $B_a$ as follow:

\begin{defi}
The stochastic embedding of the basic brick $B_a$ is given by
\begin{equation}
{\cal B}_a ={\cal D} \circ a (\bullet ,t).
\end{equation}
\end{defi}

However, classical properties of the differential calculus allow us to write $B_a$ equivalently as
\begin{equation}
B_a \cdot x= a'(x) \di {dx\over dt} .
\end{equation}
The stochastic embedding of this new form of $B_a$ is given by
\begin{equation}
\mathbb{B}_a .X=a'(X) {\cal D} X .
\end{equation}

The main problem is that in general, we do not have
\begin{equation}
\label{harmon} {\cal B}_a =\mathbb{B}_a ,
\end{equation}
as in the classical case.\\

This reflects the fact that $\stoc$ acts on operators of a given form and not on operators as an abstract element of
a given algebra. In particular, this is not a mapping.\\

Nevertheless, there exists a class of functions $a$ such that equation
(\ref{harmon}) is valid:

\begin{lem}
Equation (\ref{harmon}) is satisfied on the set $\Lambda_d$ with constant diffusion if
$a$ is an harmonic function.
\end{lem}

\begin{proof}
This follows easily from corollary \ref{corcomm}.
\end{proof}

In the sequel we study some basic properties of this embedding procedure on
differential equations.

\section{First examples}

\subsection{First order differential equations}

Let us consider a first order differential equation
$$\di {dx\over dt} =f(x,t),\eqno{1-(ODE)}$$
where $x\in \rR$ and $f :\rR \times \R \rightarrow \R$ is a given function. The
stochastic embedding of (1-ODE) leads to
$${\cal D} X=F(X,t),\eqno{1-(SODE)}$$
where $F$ is real valued.\\

The reality of $F$ imposes important constraints on solutions of 1-(SODE).
Indeed, we must have
$$DX=D_* X,$$
so that $X$ belongs to the class of Nelson-differentiable processes.\\

In our general philosophy, ordinary differential equations are only coarse
approximations to reality which must include stochastic behaviour in its
foundation. A stochastic perturbation of a first order differential equation
is then highly non-trivial. Indeed, we must consider SODE's of the form
$${\cal D} X=F(X,t)+\epsilon G(X,t),$$
where $G(X,t)$ is now complex valued. As a consequence, we allow solutions to
leave the Nelson-differentiable class.

\subsection{Second order differential equations}
\label{secondorder}

Let us consider a second order differential equation
$$\di {d^2 x\over dt^2} +a(x) \di {dx\over dt} +b(x) =0, \eqno{(2-(ODE)}$$
where $x\in \rR$, and $a,b:\rR \rightarrow \rR$ are given functions. The
stochastic embedding of $(2-(ODE))$ leads to
$${\cal D}^2 X +a(X) {\cal D} X +b(X) =0 .$$
In this case, contrary to what happens for first order differential
equations, we have no reality condition
which constrains our stochastic process.\\

In order to study such kind of equations, one can try to reduce it to a first
order equation, using standard ideas. We denote by $Y={\cal D} X$, then the
second order equation is equivalent to the following system of first order
stochastic differential equations:
\begin{equation}
\label{embedsys} \left \{
\begin{array}{lll}
{\cal D} X & = & Y,\\
{\cal D} Y & = & -a(X) Y -b(X).
\end{array}
\right .
\end{equation}
One must be careful to take $Y\in \Cu_{\C} (I)$ as $Y$ is a priori a complex
stochastic process. This remark is of importance since if we apply the
stochastic embedding procedure\footnote{Note that we have not defined the
stochastic embedding procedure on systems of differential equations.} to the
classical system of first order differential equations
\begin{equation}
\left \{
\begin{array}{lll}
\di {dx\over dt} & = & y,\\
\di {dy\over dt} & = & -a(x)y -b(x),
\end{array}
\right .
\end{equation}
by saying that we apply separately the stochastic embedding on each
differential equations, we obtain the stochastic
equation (\ref{embedsys}) but with $Y\in \Cu (I)$, which imposes strong constraints on the solutions of our equations.\\

This example proves that the stochastic embedding procedure is not so easy to
define if one wants to deal with systems of differential equations. We will
return on this problem concerning the stochastic embedding of Hamiltonian
systems.

\chapter{Reversible stochastic embedding}
\label{chapreversible}

\section{Reversible stochastic derivative}

In our construction of the stochastic derivative, we have imposed some
constraints as for example the gluing to the classical derivative on
differentiable deterministic processes. We have moreover kept some properties
of the classical derivative such as linearity. However, we have not conserved more
important properties of the classical derivative which are used in the study of
classical differential equations. For example, let us consider
$$\di {d^2 x\over dt^2} =f(x) ,\eqno{(E)}$$
which is the basic equation of Newton's mechanics. An important property of
this kind of equations is its {\it
reversibility}:\\

Let $t\rightarrow x(t)$ be a solution of (E). We denote by
$\tilde{x}(t)=x(-t)$. Then, we have
$${d^2 \tilde{x}\over dt^2} =\di{d\over dt} (-\di {dx\over dt} (-t))=\di {d^2 x\over dt^2} (-t)=f(x(-t))=f(\tilde{x}(t)),$$
proving that the reversed solution $\tilde{x} (t)$ is again a solution of the
same equation. In this case, we say that
the differential equation is {\it reversible}.\\

The reversibility argument used the following important property:
$$\di {d\over dt} (x(-t))=-\di {dx\over dt} (-t).\eqno{(R)}$$

The natural way to introduce a notion of {\it reversibility} is then to look
for the stochastic differential equation satisfied by $\tilde{X} (t)=X(-t)\in
\Cu (I)$ the reversed processes. However, in general, we do not have access to
$D\tilde{X}$ or $D_* \tilde{X}$. As a consequence, a definition using this
characterization is not effective. In the following, we follow a
different strategy.\\

A convenient way to characterize the reversibility of a given differential
equation, described by a differential operator
\begin{equation}
O=\di\sum_i a_i \di {d^i \over dt^i} \in \rR  [ d/dt ]
\end{equation}
is to prove that this operator is invariant under the substitution
\begin{equation}
r: \rR [d/dt ]  \longrightarrow \rR [d/dt ]
\end{equation}
which is $\rR$ linear and defined by
\begin{equation}
\label{fondarev} r (d/dt) =-d/dt .
\end{equation}
We then introduce in our setting, the following analogous substitution:

\begin{defi}
The reversibility operator $R :\cC [D,D_* ] \rightarrow \cC [D,D_* ]$ is a
$\cC$ morphism defined by
\begin{equation}
R(D)=-D_* ,\ \ R(D_* )=-D .
\end{equation}
\end{defi}

We have the following immediate consequence of the definition:

\begin{lem}
The reversibility operator is an involution of $\cC [D,D_* ]$.
\end{lem}

This operator acts non trivially on our stochastic derivative. Precisely, we
have:

\begin{lem}
\begin{equation}
\label{fondarevstoc}
R({\cal D} )=-\overline{\cal D} .
\end{equation}
\end{lem}

The complex nature of the stochastic derivative induces new phenomenon which
are different from the classical case. For example, we have
\begin{equation}
R({\cal D}^2 )=\overline{\cal D}^2 ,
\end{equation}
contrary to what happens for $r$.\\

We now define our notion of a reversible stochastic equation.

\begin{defi}
\label{reversi}
[Reversibility] Let $O\in \rR [D,D_* ]$, then the stochastic
equation $O\cdot X=0$ is reversible if and only if $R(O)\cdot X=0$.
\end{defi}

A natural problem is the following:\\

{\bf Reversibility problem}: {\it Find an operator such that the stochastic
embedding of a reversible
equation is again a reversible equation in the sense of definition \ref{reversi}.}\\

Let us consider the family of stochastic derivatives ${\cal D}_{\mu}$,
$\mu=0,\pm 1$. Without assuming a particular form for the underlying equation,
the preservation of the reversible character reduces to prove that the operator
$\delta$ which is chosen satisfies
\begin{equation}
R(\delta )=-\delta .
\end{equation}
In the family of stochastic derivatives ${\cal D}_{\mu}$, $\mu=0,\pm 1$, only
one case is possible:

\begin{lem}
A reversibility of a differential equation is always preserved under a
stochastic embedding if and only if this embedding is associated to the
stochastic derivative ${\cal D}_0$.
\end{lem}

\begin{proof}
Essentially this follows from equation (\ref{fondarevstoc}). If we want to
preserve reversibility then the operator ${\cal D}_{\mu}$ must satisfied
$R({\cal D}_{\mu} )=-{\cal D}_{\mu}$. This is only possible if ${\cal D}_{\mu}$
is real, i.e. $\mu =0$.
\end{proof}

It must be pointed out that the operator
$${\cal D}_0 =\di {D+D_* \over 2},$$
has been obtained by different authors using the following argument:\\

If we use only $D$ (or $D_*$) then, we give a special importance to the future
(or past) of the process, which has no physical justification. As a consequence, one
must construct an operator which combines these two quantities in a more or less
symmetric way. The simplest combination is a linear one $aD+bD_*$ with equal
coefficients $a=b$. The gluing to the classical
derivative leads to $a=b=1/2$.\\

The problem with this construction is that this argument is used on diffusion
processes, where $D$ and $D_*$ are not {\it free}. As a consequence, working
with $D$ is the same (even if the connection with $D_*$ is not trivial) than
working with $D_*$. We can not really justify then the use of ${\cal D}_0$. It
must be pointed out that E. Nelson \cite{ne1} does not use ${\cal D}_0$ in his
derivation of
the Schr\"odinger equation, but simply $D$.\\

Here, this operator is obtained by specialization of ${\cal D}_{\mu}$, which
form is imposed by our construction (linearity, gluing to the classical
derivative, reconstruction property). The reconstruction property imposes that
$\mu\not=0$ unless we
work with diffusion processes.\\

Imposing a new constraint on the reversibility on this operator leads us to
$\mu=0$. The operator ${\cal D}_0$ is of course defined on $\Cu (I)$, but in
order to satisfy the whole
constraints of our construction, we must restrict its domain to diffusion processes.\\

We can of course find reversible equations without using ${\cal D}_0$ but ${\cal D}_{\mu}$. We keep the notations and
conventions of chapter \ref{chapembed}. We first define the action of $R$ on a given operator of the form
\begin{equation}
\label{stocopera}
{\cal O}=\di\sum_{i=0}^n  a_i (\bullet ,t) (-1)^i \overline{{\cal D}}^i .
\end{equation}

\begin{defi}
The action of $R$ on (\ref{stocopera}) is denoted $R({\cal O})$ and defined by
\begin{equation}
R({\cal O}) =\di\sum_{i=0}^n a_i (\bullet ,t) {\cal D}^i .
\end{equation}
\end{defi}

The definition \ref{reversi} of a reversible equation can then be extended to cover operators of the form \ref{stocopera}.\\

Using this definition, we can prove that the stochastic equation
$$
{\cal D}_{\mu}^2 X=-\nabla U (X) ,\eqno{(E)}
$$
is reversible.\\

Indeed, we have:

\begin{lem}
Equation (E) is reversible.
\end{lem}

\begin{proof}
We have
\begin{equation}
\left .
\begin{array}{lll}
R({\cal D}_{\mu}^2 X +\nabla U (X) ) & = & \overline{\cal D}^2 X +\nabla U (X) ,\\
 & = & \overline{{\cal D}_{\mu}^2 X} +\nabla U (X) .
\end{array}
\right .
\end{equation}
As $U$ is real valued and $X$ are real stochastic processes, we deduce from (E)
that
\begin{equation}
\overline{{\cal D}_{\mu}^2 X} =-\overline{\nabla U (X)}=-\nabla U (X) .
\end{equation}
We deduce that
\begin{equation}
R({\cal D}_{\mu}^2 X +\nabla U (X) )=0,
\end{equation}
which concludes the proof.
\end{proof}

\section{Iterates}

There exists a fundamental difference between ${\cal D}_0$ and ${\cal
D}_{\mu}$, $\mu\not=0$. The operator ${\cal D}_0$ send real stochastic
processes to real stochastic processes in the contrary of ${\cal D}_{\mu}$,
$\mu\not= 0$, which leads to complex stochastic processes. As a consequence,
the $n$-i\`eme iterates of ${\cal D}_0$ is simply defined by
\begin{equation}
{\cal D}_0^n ={\cal D}_0 \circ \dots \circ {\cal D}_0 ,
\end{equation}
without problem, where a special extension of ${\cal D}_{\mu}$, $\mu\not=0$ to
complex stochastic processes must be discussed.

\section{Reversible stochastic embedding}

Using ${\cal D}_0$, we can define a stochastic embedding which conserves the
fundamental property of reversibility of a given equation. We keep notations from chapter \ref{chapembed}.

\begin{defi}[Reversible stochastization]
\label{stocarev}
Let $O\in \kk [d/dt ]$ be a differential operator, of the form
$$O=a_0 (\bullet,t)+a_1 (\bullet,t)\di {d\over dt} +\dots + a_n (\bullet,t)\di {d^n \over dt^n }, \ a_i \in \kk,\ =0,\dots ,n,$$
where $n\in \nN$ is given.

The reversible stochastic embedding of $O$ is an element $O_{\rm rev}$ of $\Cu
(I) [{\cal D}_0 ]$ defined by
\begin{equation}
O_{\rm rev} =a_0 (\bullet ,t) +a_1 (\bullet ,t)\di {\cal D}_0 +\dots + a_n (\bullet ,t) \di {{\cal D}_0^n }, \ a_i
\in \Cu (I),\ i=0,\dots ,n.
\end{equation}
\end{defi}

A differential equation (E) is defined by a differential operator $O\in \kk
[d/dt ]$, i.e. an equation of the form
$$O\cdot x=0 , \eqno{(E)}$$
where $x$ is a function.

Using stochastization, the reversible stochastic analogue of (E) is defined by
$$O_{\rm rev} \cdot X=0, \eqno{(RSE)}$$
where $X$ is a stochastic process.

\section{Reversible versus general stochastic embedding}

The reversible stochastic embedding leads to very different results than the
general stochastic embedding. We can already see this difference on first order
differential equations. Let us consider
$$\di {dx\over dt}=f(x),$$
where $x\in \rR$ and $f$ is a real valued function. The reversible stochastic
embedding gives
$${\cal D}_0 X=f(X).$$
Contrary to what happens for the stochastic embedding, this equation does not
impose for the solution to be a Nelson differentiable processes.

\section{Stochastic mechanics and the Stochastization procedure}

\subsection{The Stochastic Newton Equation}

The stochastized version of the classical system:
\begin{eqnarray}
\dot{x}(t) & = & v(t)\nonumber\\
\dot{v}(t) & = & K(x(t))
\end{eqnarray}
is given by:
\begin{eqnarray}\label{NSE}
    \mathcal{D}X(t) & = & V(t)\nonumber\\
\mathcal{D}V(t) & = & K(X(t))\label{ipart}
\end{eqnarray}
where $V\in\Cu_{\C} (I)$ and $K$ is a force: $K(x)=-\nabla U(x)$ and $U$ a potential.\\

We can give at least two different kind of solutions of this equation, and so two relevant models.\\

In the first one, the component $X$ is the position in the Ornstein-Uhlenbeck
theory of Brownian Motion and is not submitted to a random noise. The system
writes:
\begin{equation}
    \left\{
\begin{array}{l}
    dX(t)=V(t)dt\\
    dV(t)=-\alpha V(t)dt+K(X(t))dt+\sigma dW(t) \\
    X(0)=X_0,\ V(0)=V_0,
\end{array}\right.
\end{equation}

We have noticed in a previous section that, at an equilibrium (i.e. $X$ has a
stationary density) and if $e^{-U}$ is integrable, then:
\begin{eqnarray}
    \mathcal{D}X(t) & = & V(t),\\
    \mathcal{D}^2 X(t) & = & K(X(t))-i\alpha V(t).
\end{eqnarray}
Therefore $(X,V)$ solves the Newton stochastized system (\ref{NSE}) if and only
if $\alpha=0$. Moreover we note in this
particuliar case that $X$ is a Nelson-differentiable process.\\

The second one is described by
\begin{equation}
dX(t)=b(t,X(t))dt+\sigma dW(t) ,
\end{equation}
where the function $b$ must be determined. In this case, we proved that the density
$p_t(x)$ of a solution $X$ of (\ref{NSE}) writes
$p_t(x)=\Psi(t,x)\overline{\Psi}(t,x)$ where $\Psi$ solves the Schr\"odinger
equation: $i\sigma^2\partial_t\Psi+\frac{\sigma^4}{2}\partial_{xx}\Psi=U\Psi$.
In this case, $X$ is driven by a Brownian motion and is not
Nelson-differentiable.

\part{Stochastic embedding of Lagrangian and Hamiltonian systems}

\chapter{Stochastic Lagrangian systems}
\setcounter{section}{0} \setcounter{equation}{0}

Most of classical mechanics can be formulated using Lagrangian formalism
(\cite{ar},\cite{am}). Lagrangian mechanics contains important
problems, like the $n$-body problem. Using our framework, we study Lagrangian
dynamical systems under stochastic perturbations\footnote{For the $n$-body
problem, which is usually used to study the long term behavior of the solar
system \cite{marmi}, this problem is of crucial importance. Indeed, the
$n$-body problem is only an approximation of the real problem, and even if some
numerical simulations take into account relativistic effects \cite{laskar}, this is not sufficient \cite{moser}.}.

Our approach is first to embed classical Lagrangian systems, in particular the
associated Euler-Lagrange equation (EL) in order to obtain an idea of what kind
of equation govern stochastic Lagrangian systems. We then develop a stochastic
calculus of variations. We obtain an analogue of the {\it least-action
principle}\footnote{In our case, the word {\it least-action} is misleading and
a better terminology is {\it stationary} (see below).} which gives a second
stochastic Euler-Lagrange equation, denoted by (SEL) in the sequel. We then
prove the following surprising result, called the coherence lemma: we have
$\stoc (EL)=(SEL)$.

The principal interest of Lagrangian systems is that the action of a group of
symmetries leads to first integrals of motion, i.e. functions which are constants
on solutions of the equations of motion. The celebrated theorem of E. Noether
gives a precise relation between symmetries and first integrals. We prove a
stochastic analogue of E. Nother theorem.

Finally, we prove that the stochastic embedding of Newton's Lagrangian systems
lead to a non linear Schr\"odinger's equation for a given wave function whose
modulus is equal to the probability density of the underlying stochastic
process.

\section{Reminder about Lagrangian systems}

We refer to \cite{ar} for more details, as well as \cite{am}.\\

Lagrangian systems play a central role in dynamical systems and physics, in
particular for mechanical systems. A Lagrangian system is defined by a
Lagrangian function, commonly denoted by $L$, and depending on three variables:
$x$, $v$, and $t$ which belongs in the sequel to $\rR$. As Lagrangian systems
come from mechanics, the letter $x$ stands for position, the letter $v$ for
speed and the letter $t$ for time. In what follows, we consider a special type
of Lagrangian function called admissible in the following.

\begin{defi}
An admissible Lagrangian function is a function $L$ such that:\\

{\rm i)} The function $L(x,v,t)$ is defined on $\rR^d \times \cC^d \times \rR$,
holomorphic in the second variable and real for $v\in \rR$.

{\rm ii)} $L$ is autonomous, i.e. $L$ does not depend on time.
\end{defi}

Condition i) is fundamental. This condition is necessary in order to apply the
stochastization procedure (see below). The fact that we only consider
autonomous Lagrangian function is due to technical difficulties in order to
take into account backward and forward filtrations in the computation of the
stochastic Euler-Lagrange equation (see below).

\begin{rema}
In applications, admissible Lagrangian functions $L$ are analytic extensions to
the complex domain of real analytic Lagrangian functions. For example, the
classical Newtonian Lagrangian $L(x,v)=(1/2)v^2 -U(x)$, defined on an
open\footnote{This Lagrangian function is not always defined on $\R \times \R$.
An example is given by Newton's potential $U(x)=1/x$, $x\in\R^*$.} subset of
$\R \times \R$, with an analytic potential is an admissible Lagrangian
function.
\end{rema}

A Lagrangian function $L$ being given, the equation
$$
\di {d\over dt} \left ( \di {\partial L \over \partial v} \right ) =\di
{\partial L\over \partial x} . \eqno{(EL)}
$$
is called the {\it Euler-Lagrange equations}.\\

An important property of the Euler-Lagrange equation is that it derives from a {\it
variational principle}, namely the {\it least action principle} (see
\cite{ar},p.59). Precisely, a curve $\gamma : t\mapsto x(t)$ is an {\it
extremal}\footnote{We refer to \cite{ar}, chapter 3, $\S$.12 for an
introduction to the {\it calculus of variations}.} of the functional
$$J_{a,b} (\gamma)=\di\int_a^b L(x(t),\dot{x}(t) ,t) dt ,$$
on the space of curves passing through the points $x(a)=x_a$ and $x(b)=x_b$, if
and only if it satisfies the Euler-Lagrange equation along the curve $x(t)$.

\section{Stochastic Euler-Lagrange equations}

We now apply our stochastic procedure $\stoc$ to an admissible Lagrangian.

\begin{lem}
Let $L(x,v) :\rR^d \times \cC^d \rightarrow \cC$ be an admissible Lagrangian
function. The stochastic Euler-Lagrange equation obtained from (EL) by the
stochastic procedure is given by
$$
\di {\cal D}_{\mu} \di\left (  \di {\partial L\over \partial v} (X(t),{\cal
D}_{\mu} X(t) \right ) =\di {\partial L \over \partial x} (X(t),{\cal D}_{\mu}
X(t)) . \eqno{\stoc(EL)}
$$
\end{lem}

\begin{proof}
The Euler-Lagrange equation associated to $L(x,v)$ can be seen as the following
differential operator
$$O_{EL} =\di {d\over dt}\circ \di {\partial L\over \partial v} -\di {\partial L\over \partial x} ,$$
acting on $(x(t),\dot{x}(t))$. The embedding of $O_{EL}$ gives
$${\cal O}_{EL}=\di {\cal D}_{\mu} \circ \di {\partial L\over \partial v} -\di {\partial L\over \partial x} .$$
As $O_{EL}$ acts on $(x(t),\dot{x}(t))$, the operator ${\cal O}_{EL}$ acts on
$(X(t), {\cal D}_{\mu} X (t))$. This concludes the proof.
\end{proof}

The free parameter $\mu\in \{ -1,0,1\}$ can be fixed depending on the nature of the extension used.\\

It must be pointed out that there exist crucial differences between all these
extensions due to the fact that ${\cal D}_{\mu}$ is complex valued for $\mu=\pm
1$ and real for $\mu=0$. Indeed, let us consider the following admissible
Lagrangian function:
$$L(x,v)=\di {1\over 2} v^2 -U(x),$$
where $U$ is a smooth real valued function. Then, equation $\stoc$(EL) gives
$${\cal D}_{\mu} V= U(X),$$
where $V={\cal D}_{\mu} X$. When $\mu=\pm 1$, this equation imposes strong
constraints on $X$ due to the real
nature of $U(X)$, namely that ${\cal D}^2_{\mu} X\in \Nd (I)$.\\

On the contrary, when $\mu=0$, i.e. in the {\it reversible case}, these {\it
intrinsic} conditions disappear.

\section{The coherence problem}

Up to now, the stochastic embedding procedure can be viewed as a formal manipulation of differential equations. Moreover, as
most classical manipulations on equations do not commute with the stochastic embedding, this procedure is not canonical
\footnote{We return to this problem in our discussion of a stochastic symplectic geometry which can be used to bypass this
kind of problem.}. In order to rigidify this construction and to make precise the role of this stochastic embedding procedure,
we study the following problem, called {\it the coherence problem}:\\

We know that the Euler-Lagrange equations are obtained via a least-action principle on a functional. The main problem
is the existence of a stochastic analogue of this least-action principle, that we can call a {\it stochastic least action
principle}, compatible with the stochastic embedding procedure.

\begin{eqnarray}
\xymatrix{
  & L(x(t),\dot{x}(t)) \ar[d]_{\mbox{\rm Least action principle}} \ar[r]^{\stoc} & L(X(t),{\cal D}X(t))
  \ar[d]^{\mbox{\rm Stochastic least action principle ? }}       \\
  & (\mbox{\rm EL}) \ar[r]_{\stoc}   & (\mbox{\rm SEL})               }
\end{eqnarray}

In the next chapter, we develop the necessary tools to answer to this problem, i.e. a {\it stochastic calculus of
variations}. Note that due to the fact that the stochastic Lagrangian as well as the stochastic Euler-Lagrange equation
are fixed, this problem is far from being trivial. The main result of the next chapter is the Lagrangian coherence lemma which
says precisely that the stochastic Euler-Lagrange equation obtained via the stochastic embedding procedure coincide with the
characterization of extremals for the functional associated to the stochastic Lagrangian function using the stochastic
calculus of variations. As a consequence, we obtain a rigid picture involving the stochastic embedding procedure and a first
principle via the stochastic least action principle.\\

This picture will be then extended in another chapter when dealing with the Hamiltonian part of this theory.

\chapter{Stochastic calculus of variations}
\label{chapvar}

The embedding procedure allows us to associate a stochastic Euler-Lagrange
equation to a stochastic Lagrangian function. A basic question is then the
existence of an analogue of the {\it least action principle}. In this section,
we develop a stochastic calculus of variations for our Lagrangian function
following a previous work of K. Yasue \cite{ya1}. Our main result, called the
coherence lemma, states that the stochastic Euler-Lagrange equation can be obtained
as an application of a stochastic least action principle. Moreover, this derivation is
consistent with the stochastic embedding procedure.

\section{Functional and $L$-adapted process}

In the sequel we denote by $I$ a given open interval $(a,b)$, $a<b$.\\

We first define the stochastic analogue of the classical functional.

\begin{defi}
Let $L$ be an admissible Lagrangian function. The functional associated to $L$
is defined by
\begin{eqnarray}
\label{functional} J_{a,b}(X)=E\left[\int_a^b L(X(t),\mathcal{D}_{\mu}
X(t))dt\right] ,
\end{eqnarray}
for all $X\in \Cu (I)$.
\end{defi}

In what follows, we need a special notion introduced by Yasue \cite{ya1}, and
called $L$-adaptation:

\begin{defi}
Let $X\in \Cu (I)$ be a stochastic process. We denote by $\cal P$ and $\cal F$
the past and the future of $X$. Let $L$ be an
admissible Lagrangian function. A process $X\in \Cu (I)$ is called $L$-adapted if:\\

{\rm i)} $\di {\partial L\over \partial v} (X(t),{\cal D}_{\mu} X(t))$ is adapted to
$\cal P$ and $\cal F$.

{\rm ii)} $\di {\partial L\over \partial v} (X(t),{\cal D}_{\mu} X(t)) \in\Cu (I)$.
\end{defi}

Diffusion processes are $L$-adapted.\\

\section{Space of variations}

Calculus of variations is concerned with the behaviour of functionals under
{\it variations} of the underlying functional space, i.e. objects of the form
$\gamma +h$, where $\gamma$ belongs to the functional space and $h$ is a given
functional space of variations. A special care must be taken in our case to
define what is the class of {\it variations} we are considering. In general,
this problem is not really pointed out as both variations and curves can be
taken in the same functional space (see \cite{ar},p.56,footnote 26). We
introduce the following terminology:

\begin{defi}
Let $P$ be a subspace of $\Cu (I)$ and $X\in \Cu (I)$. A $P$-variation of $X$
is a stochastic process of the form $X+Z$, where $Z\in P$.
\end{defi}

In the sequel, we consider two subspaces of variations: $\Nd (I)$ and $\Cu (I)$.\\

The choice of $\Cu (I)$ is natural. However, doing this we can obtain
stochastic processes with completely different behaviour than $X$\footnote{Of
course, this is not the case in the classical case: one consider $x\in
C^{\infty} (I)$ and $z\in C^{\infty} (I)$ such that $x+h\in C^{\infty} (I)$ is
very similar to $x$. For example, we don't choose $z\in C^0 (I)$
which leads to radically new behaviour of $x+z$ with respect to $x$.}.\\

What is the specific property of $X\in \Cu (I)$ that we want to keep ?\\

If we refer to the construction of the stochastic derivative, then a main point
is the existence of an {\it imaginary part} in ${\cal D}_{\mu} X$\footnote{Of
course, as long as $\mu=\pm 1$. This is of importance since we will be able to
choose a more general variations space in this case.}. This property is related
to the non-differentiability of the underlying stochastic process. We are then
lead to search for variations $Z$ which conserve this imaginary part. As a
consequence, we must consider Nelson difference processes introduced in the
previous part\footnote{An analogous problem is considered in \cite{cr2}, where
a non differentiable variational principle is defined.}, and denoted by $\Nd
(I)$.

\section{Differentiable functional and stationary processes}

We now define our notion of {\it differentiable functional}. Let $P$ be a
subspace of $\Cu (I)$.

\begin{defi}
Let $L$ be an admissible Lagrangian function and $J_{a,b}$ the associated
functional. The functional $J_{a,b}$ is called $P$-differentiable at an
$L$-adapted process $X\in\Cu (I)$ if
\begin{equation}
J_{a,b} (X+Z)-J_{a,b} (X)=dJ_{a,b} (X,Z)+R(X,Z) ,
\end{equation}
where $dJ_{a,b}(X,Z)$ is a linear functional of $Z\in P$ and
$R(X,Z)=o(\parallel Z\parallel )$.
\end{defi}

The stochastic analogue of a stationary point is then defined by:

\begin{defi}
A $P$-stationary process for the functional $J_{a,b}$ is a stochastic process
$X\in\Cu (I)$ such that $dJ(X,Z)=0$ for all $Z\in P$.
\end{defi}

\subsection{The $P=\Cu (I)$ case}

Our main result is:

\begin{lem}
The functional $J_{ab}$ defined by (\ref{functional}) is $\Cu
(I)$-differentiable at any $L$-adapted process $X\in\Cu (I)$, and for
all $Z\in\Cu (I)$, the differential is given by:
\begin{equation}
\label{differentielle} \left .
\begin{array}{lll}
dJ_{ab}(X,Z) & = & E\left [ \di\int_a^b\left [ \di {\partial L\over \partial
x}(X(u),\mathcal{D}_{\mu}X(u))- \mathcal{D}_{-\mu}\left ( \di {\partial L\over
\partial v}
(X(u),\mathcal{D}_{\mu}X(u)) \right ) \right ] Z(u)du \right ] \\
 & & +g(Z,\partial_v L)(b)-g(Z,\partial_v L)(a) ,
\end{array}
\right .
\end{equation}
where
\begin{equation}
g(Z,\partial_v L)(s)=E\left[Z(u)\partial_v L(X(u),\mathcal{D}_{\mu}X(u))\right]
.
\end{equation}
\end{lem}

\begin{proof}
Let $X$ and $Z$ be two $L$-adapted processes. The Taylor expansion of $L$
gives:
\begin{equation}
\left .
\begin{array}{lll}
L(X+Z,\mathcal{D}_{\mu}(X+Z))-L(X,\mathcal{D}_{\mu}(X)) & = & \partial_x L(X,\mathcal{D}_{\mu}(X))Z \\
 & & +\partial_v L(X,\mathcal{D}_{\mu}(X))\mathcal{D}_{\mu}(Z)\\
 & & +o(\left\|Z\right\|) ,
\end{array}
\right .
\end{equation}
which yields (\ref{differentielle}) by integration and (\ref{formule2}).
\end{proof}

\subsection{The $P=\Nd (I)$ case}

Our main result is:

\begin{lem}
The functional $J_{ab}$ defined by (\ref{functional}) is $\Nd
(I)$-differentiable at any $L$-adapted process $X\in \Cu (I)$, and for
all $Z\in \Nd (I)$ the differential is given by:
\begin{equation}
\label{differentielle} \left .
\begin{array}{lll}
dJ_{ab}(X,Z) & = & E\left[\di \int_a^b(\partial_x L-\mathcal{D}_{\mu}\partial_v L)(X(u),\mathcal{D}_{\mu}X(u))Z(u)du\right] \\
 & & +g(Z,\partial_v L)(b)-g(Z,\partial_v L)(a) ,
\end{array}
\right .
\end{equation}
where
\begin{equation}
g(Z,\partial_v L)(s)=E\left[Z(u)\partial_v L(X(u),\mathcal{D}_{\mu}
X(u))\right] .
\end{equation}
\end{lem}

\begin{proof}
Let $X$ a $L-$adapted process and $H$ a Nelson-differentiable process. The
Taylor expansion of $L$ gives
\begin{equation}
\left .
\begin{array}{lll}
L(X+H,\mathcal{D}_{\mu}(X+H))-L(X,\mathcal{D}_{\mu}(X)) & = & \partial_x L(X,\mathcal{D}_{\mu}(X))H\\
& & +\partial_V L(X,\mathcal{D}_{\mu}(X))\mathcal{D}_{\mu}(H)\\
 & & +o(\left\|H\right\|) ,
\end{array}
\right .
\end{equation}
which yields (\ref{differentielle}) by integration and (\ref{derivation}).
\end{proof}

\section{A technical lemma}

The classical derivation of the least action principle used a well known result about {\it bump functions} (see
\cite{ar},p.57). In the stochastic framework, we will need the following result:

\begin{lem}
\label{plateau}
Let $Y\in {\cal P}_{\C}$ be a complex stochastic process. If $Y$ satisfies
\begin{equation}
\label{equaplateau}
\int_0^1 E\left [ Y(u) {\cal D}_{\mu} Z (u) \right ] \, du =0 ,
\end{equation}
for all $Z\in {\cal N}^1 ([0,1])$ then $Y$ is a constant process.
\end{lem}

\begin{proof}
We denote $Y=Y_1 +iY_2$, where $Y_i \in {\cal P}_{\R}$ and ${\cal D}_{\mu} Z =A$, where $A\in {\cal P}_{\R}$. The
equation (\ref{equaplateau}) is equivalent to
\begin{equation}
\label{equaplateau2}
\left .
\begin{array}{lll}
\int_0^1 E \left [ Y_1 (u)A(u)\right ] \, du & = & 0 ,\\
\int_0^1 E \left [ Y_2 (u)A(u)\right ] \, du & = & 0,
\end{array}
\right .
\end{equation}
for all $A\in {\cal P}_{\R}$ such that there exists $Z\in \Cu ([0,1])$ satisfying ${\cal D}_{\mu} Z =A$.\\

Let $Z_{Y_1}$ be the process defined by
\begin{equation}
Z_{Y_1} (u)=\int_0^u Y_1 (s)ds -u\int_0^1 Y_1 (s)ds .
\end{equation}
We have $Z_{Y_1} \in {\cal N}^1 (I)$ with $Z(0)=Z(1)=0$. Indeed, we have
\begin{equation}
{\cal D}_{\mu} Z (u)= Y_1 (u) -\int_0^1 Y_1 (s)\, ds .
\end{equation}
As a consequence, we have in our notations $B=0$ and the first equation of (\ref{equaplateau2}) reduces to
\begin{equation}
\int_0^1 E[Y_1(u) A(u)]du = E\left[\int_0^1\left(Y_1 (u)-\int_0^1 Y_1 (s)ds\right)^2du\right].
\end{equation}
We deduce that $Y_1$ is a constant process, that is for all $u\in I$, $Y_1 (u)=C$ a.s., where $C$ is a random variable.\\

The same argument with the second equation of (\ref{equaplateau2}) and $Z_{Y_2}$ concludes the proof of the lemma.
\end{proof}

\section{Least action principles}

As for the computation of the differential of functionals, we must consider two
cases: $P=\Cu (I)$ and $P=\Nd (I)$.

\subsection{The $P=\Cu (I)$ case}

The main result of this section is the following analogue of the least-action
principle for Lagrangian mechanics.

\begin{thm}[Global Least action principle]
A necessary and sufficient condition for an $L$-adapted process to be a $\Cu
(I)$-stationary process of the functional $J_{ab}$ with fixed end points
$X(a):=X_a\in H$ et $X(b):=X_b\in H$ is that it satisfies
\begin{eqnarray}
{\partial L\over \partial x}(X(t),\mathcal{D}_{\mu} X(t))-\mathcal{D}_{-\mu}
\left [ \di {\partial L\over \partial v} (X(t),\mathcal{D}_{\mu} X(t)) \right ]
=0 .
\end{eqnarray}
We call this equation the Global Stochastic Euler-Lagrange equation (GSEL).
\end{thm}

We have conserved the terminology of least-action principle even if we have no
notion of {\it extremals} for our complex valued functional.

\begin{proof}
We denote by $I=]0,1[$. Let $X\in\EuScript{C}^1(I)$ be a solution of
\begin{equation}
\d (\partial_x \mathcal{L}-\mathcal{D}_\mu\partial_{v} L)(X(u),\mathcal{D}_\mu
X(u))=0 ,
\end{equation}
then $X$ is a $\mathcal{N}^1(I)$-stationary process for the functional $J_I$.\\

Conversely, let $X$ is a $\Cu (I)$-stationary process for the functional $J_I$,
i.e. $dJ_I(X,Z)=0$. Writing
$$ (\partial_x
\mathcal{L}-\mathcal{D}_\mu\partial_{v} L)(X(u),\mathcal{D}_\mu
X(u))=\mathcal{D}_\mu Y(u), $$
where
\begin{equation}
Y(u)=\int_0^u \partial_x L(X(s),\mathcal{D}_\mu X(s))ds -\partial_{v}
L(X(u),\mathcal{D}_\mu X(u)),
\end{equation}
we obtain for any $Z\in \Cu (I)$ with $Z(0)=Z(1)=0$:
\begin{eqnarray}
dJ_I(X,Z) & = & E\left[\int_0^1 \mathcal{D}_\mu Y(u) Z(u)du\right]\\ \nonumber
          & = & \int_0^1 E[\mathcal{D}_\mu Y(u) Z(u)]du .\\ \nonumber
\end{eqnarray}
Using the $\Cu (I)$-product rule (see equation \ref{formule2}), we obtain
\begin{equation}
dJ_I(X,Z)=-\int_0^1 E[Y(u) \mathcal{D}_\mu Z(u)]du .
\end{equation}

Using lemma \ref{plateau} we obtain that $Y$ is a constant process.\\

Hence, we have $\mathcal{D}_\mu Y(u)=0$ and
\begin{equation}
(\partial_x \mathcal{L}-\mathcal{D}_\mu\partial_{v} L)(X(u), \mathcal{D}_\mu
X(u))=0 ,
\end{equation}
which concludes the proof.
\end{proof}

\subsection{The $P=\Nd (I)$ case}

Our main result is:

\begin{thm}[least action principle]
A necessary and sufficient condition for an $L$-adapted process to be a $\Nd
(I)$-stationary process of the functional $J_{ab}$ with fixed end points
$X(a):=X_a\in H$ et $X(b):=X_b\in H$ is that it satisfies
\begin{eqnarray}
{\partial L\over \partial x}(X(t),\mathcal{D}_{\mu} X(t))-\mathcal{D}_{\mu}
\left [ \di {\partial L\over \partial v} (X(t),\mathcal{D}_{\mu} X(t)) \right ]
=0 .
\end{eqnarray}
We call this equation the weak stochastic Euler-Lagrange equation (SEL).
\end{thm}

\begin{proof}
We denote by $I=]0,1[$. Let $X\in\EuScript{C}^1(I)$ be a solution of
\begin{equation}
\d (\partial_x \mathcal{L}-\mathcal{D}_\mu\partial_{v} L)(X(u),\mathcal{D}_\mu
X(u))=0 ,
\end{equation}
then $X$ is a $\mathcal{N}^1(I)$-stationary process for the functional $J_I$.\\

Conversely, let $X$ is a $\mathcal{N}^1(I)$-stationary process for the
functional $J_I$, i.e. $dJ_I(X,Z)=0$. Writing $$ (\partial_x
\mathcal{L}-\mathcal{D}_\mu\partial_{v} L)(X(u),\mathcal{D}_\mu
X(u))=\mathcal{D}_\mu Y(u), $$ where
\begin{equation}
Y(u)=\int_0^u \partial_x L(X(s),\mathcal{D}_\mu X(s))ds -\partial_{v}
L(X(u),\mathcal{D}_\mu X(u)),
\end{equation}
we obtain for any $Z\in\mathcal{N}^1(I)$ with $Z(0)=Z(1)=0$:
\begin{eqnarray}
dJ_I(X,Z) & = & E\left[\int_0^1 \mathcal{D}_\mu Y(u) Z(u)du\right]\\ \nonumber
          & = & \int_0^1 E[\mathcal{D}_\mu Y(u) Z(u)]du .\\ \nonumber
\end{eqnarray}
Using the ${\cal N}^1 (I)$-product rule (see equation \ref{derivation}), we
obtain
\begin{equation}
dJ_I(X,Z) =-\int_0^1 E[Y(u) \mathcal{D}_\mu Z(u)]du .
\end{equation}
Using lemma \ref{plateau}, we deduce that $Y$ is a constant process, that is for all $u\in I$, $Y(u)=C$ a.s.
where $C$ is a random variable.\\

Hence, we obtain $\mathcal{D}_\mu Y(u)=0$ and
\begin{equation}
(\partial_x \mathcal{L}-\mathcal{D}_\mu\partial_{v} L)(X(u), \mathcal{D}_\mu
X(u))=0 ,
\end{equation}
which concludes the proof.
\end{proof}

\section{The coherence lemma}

It is not clear that the stochastic Euler-lagrange equation obtained by the
stochastization procedure and the $\Nd (I)$ or $\Cu (I)$ least-action principle
coincide. One easily sees that this is not the case for $P=\Cu (I)$. In the
contrary, we have the following lemma, called the {\it coherence lemma}, which
ensure that for $P=\Nd (I)$ we obtain the same equations.

\begin{lem}[coherence lemma]
The following diagram commutes :
\begin{eqnarray}
\xymatrix{
  & L(x(t),x'(t)) \ar[d]_{\mbox{\rm Least action principle}} \ar[r]^{\stoc} & L(X(t),\mathcal{D}X(t))
  \ar[d]^{\mbox{\rm Stochastic Least action principle}}       \\
  & (EL) \ar[r]_{\stoc}   & (SEL)
  }
\end{eqnarray}
\end{lem}

\begin{proof}
This is an immediate consequence of the previous results.
\end{proof}

\begin{rema}
When $\mu=0$, i.e. in the reversible case, the previous lemmas and theorems are
true under $\Cu (I)$ variations. Note that when $\mu=0$, our stochastic
derivatives coincides with the Misawa-Yasue \cite{my} canonical formalism for
stochastic mechanics.
\end{rema}

\chapter{The Stochastic Noether theorem}
\label{noetherchap}

A natural question arising from the stochastization procedure of classical
dynamical systems, in particular, Lagrangian systems, is to understand what
remains from classical first integrals of motion. First integrals play a
central role in many problems like the $n$-body problem. In this section, we
obtain a stochastic analogue of the Noether theorem. We then defined the notion
of first integrals for stochastic dynamical systems. We also discuss the
consequences of the existence of first integrals in the context of chaotic
dynamical systems.

\section{Tangent vector to a stochastic process}

Let $X\in \Cu (I)$ be a stochastic process. We define the analogue of a tangent
vector to $X$ at point $t$.

\begin{defi}
Let $X\in \Cu (I)$, $I\subset \R$. The tangent vector to $X$ at point $t$ is
the random variable ${\cal D}X (t)$.
\end{defi}

\begin{rema}
Of course, in order to define stochastic Lagrangian systems in an intrinsic
way, one must define the stochastic analogue of the tangent bundle to a smooth
manifold. In our case, it is not clear what is the adequate geometric object
underlying stochastic Lagrangian dynamics. For example, we can think of
multidimensional Brownian surfaces (\cite{fal},$\S$.16.4). All these questions
will be developed in a forthcoming paper \cite{cd2}.
\end{rema}

\section{Canonical tangent map}

In the sequel, we will need the following mapping called the {\it canonical
tangent map}:

\begin{defi}
For all $X\in \Cu (I)$, we define the canonical tangent map as
\begin{equation}
T:\ \left .
\begin{array}{lll}
\Cu (I) & \longrightarrow & \Cu (I) \times {\cal P}_{\C} ,\\
X & \longmapsto & (X,{\cal D} X ) .
\end{array}
\right .
\end{equation}
\end{defi}

The mapping $T$ will be used in the following section to define the analogue of
the linear tangent map for a stochastic suspension of a one parameter group of
diffeomorphisms.

\section{Stochastic suspension of one parameter family of diffeomorphisms}

We begin by introducing a useful notion of {\it stochastic suspension} of a
diffeomorphism.

\begin{defi}
Let $\phi :\R^n \rightarrow \R^n$ be a diffeomorphism. The stochastic
suspension of $\phi$ is the mapping $\Phi :{\cal P}\rightarrow {\cal P}$
defined by
\begin{equation}
\forall X\in {\cal P} , \ \Phi(X)_t (\omega )=\phi (X_t (\omega )) .
\end{equation}
\end{defi}

In what follows, we will frequently use the same notation for the suspension
of a given diffeomorphism and the diffeomorphism.

\begin{rema}
It seems strange that we have not defined directly the notion of diffeomorphism
on a subset $E$ of the stochastic processes, i.e. mapping $\Phi :E\rightarrow
E$ which are Fr\'echet differentiable with an inverse which is also Fr\'echet
differentiable. However, these objects do not always exist.
\end{rema}

Using the stochastic suspension, we are able to define the notion of stochastic
suspension for a one-parameter group of diffeomorphisms.

\begin{defi}
\label{group1para}
A one-parameter group of transformations $\Phi_s :
A\rightarrow A$, $s\in \R$, where $A \subset {\cal P}$, is called a
$\phi$-suspension group acting on $A$ if there exist a one parameter group of
diffeomorphisms $\phi_s :\R^n \rightarrow
\R^n$, $s\in \R$, such that for all $s\in \R$, we have:\\

i) $\Phi_s$ is the stochastic suspension of $\phi_s$,

ii) for all $X\in A$, $\Phi_s (X)\in A$.
\end{defi}

This notion of suspension group comes from our framework. It relies on the fact
that we want to understand how symmetries of the underlying Lagrangian systems
are transported via the stochastic embedding. The non-trivial condition on
the stochastic suspension of a one-parameter group of diffeomorphisms acting on $A$
comes from condition ii). However, imposing some conditions on the underlying
one parameter group, we can obtain a stochastic one parameter group which acts
on the set $E$ of good
diffusion processes.\\

Precisely, let us introduce the following class of one-parameter groups:\\

\begin{lem}
An admissible one parameter group of diffeomorphisms $\Phi=\{\phi_s\}_{s\in\R}$
is a one parameter group of $C^2$-diffeomorphisms on $\R^n$ such that
\begin{equation}
(s,x)\mapsto \frac{\partial}{\partial x}\phi_s(x)\ \mbox{\rm is}\ C^2 .
\end{equation}
\end{lem}

The main property of admissible one parameter groups is the fact that they are
well-behaved on the set of good diffusions.

\begin{lem}
\label{tecnicnoether} Let $\Phi =(\phi_s )_{s\in \R}$ be a stochastic
suspension of an admissible one parameter group of diffeomorphisms. Then,
for all $X\in E$, we have for all $t\in I$, and all $s\in \R$:\\

i) The mapping $s\mapsto D_{\mu} \Phi_s X(t) \in C^1 (\R )$, (a.s.),\\

ii) We have $\di\frac{\partial}{\partial s}[
\di\mathcal{D}_{\mu}(\phi_s(X))]=\di\mathcal{D}_{\mu} \left[ \di\frac{\partial
\phi_s(X)}{\partial s}\right]\quad (a.s.).$
\end{lem}

This lemma is trivial in the classical case where $X$ is a smooth function and
${\cal D}_{\mu}$ is the classical derivative with respect to time. Indeed, it
reduces to the Schwarz lemma. However, this inequality plays an essential role
in the derivation of the classical Noether's theorem (see \cite{ar},p.89).

\begin{proof}
According to (\ref{darsien}),
$$\mathcal{D}_\mu \phi_s(X)(t)=\mathcal{D}_\mu X(t)\cdot\frac{\partial_x \phi_s}{\partial x}X(t) +i\mu \frac{\sigma(t,X_t)^2}{2}\frac{\partial^2_x \phi_s}{\partial x^2}X(t)\quad (a.s.).$$
So : $$ \frac{\partial}{\partial s}\mathcal{D}_\mu \phi_s(X)(t)=\mathcal{D}_\mu
X(t)\cdot\frac{\partial}{\partial s}\frac{\partial_x \phi_s}{\partial x}X(t)
+i\mu \frac{\sigma(t,X_t)^2}{2}\frac{\partial}{\partial s}\frac{\partial^2_x
\phi_s}{\partial x^2}X(t)\quad (a.s.).$$ Since $(s,x)\mapsto \phi_s(x)$ is
$C^2$, we have $\d \frac{\partial}{\partial s}\frac{\partial}{\partial
x}\phi_s(x)=\frac{\partial}{\partial x}\frac{\partial}{\partial s}\phi_s(x)$ by
the Schwarz lemma.
In the same way, $$ \frac{\partial}{\partial s}\frac{\partial^2}{\partial x^2}\phi_s(x)=\frac{\partial}{\partial x}\frac{\partial}{\partial s}\frac{\partial}{\partial x}\phi_s(x)$$ because $(s,x)\mapsto \frac{\partial}{\partial x}\phi_s(x)$ is $C^2$.\\
Therefore  : $$ \frac{\partial}{\partial s}\frac{\partial^2}{\partial
x^2}\phi_s(x)=\frac{\partial^2}{\partial x^2}\frac{\partial}{\partial
s}\phi_s(x).$$
Applying (\ref{darsien}) to $\d \frac{\partial}{\partial
s}\phi_s$, we can conclude that :
$$\frac{\partial}{\partial s}[\mathcal{D}_{\mu}(\phi_s(X))]=\mathcal{D}_{\mu}\left[\frac{\partial \phi_s(X)}{\partial s}\right]\quad (a.s.).$$
\end{proof}

It must be pointed out that every extension of this lemma will lead to a
substantial improvement of the following stochastic Noether theorem.

\section{Linear tangent map}

Let $X\in \Cu (I)$ and $\phi : \R^n  \rightarrow \R^n$ be a diffeomorphism. The
image of $X$ under the stochastic suspension of $\phi$, denoted by $\Phi$,
induces a natural map for tangent vectors denoted by $\Phi_*$, called the {\it
linear tangent map}, and defined as in classical differential geometry by:

\begin{defi}
Let $\Phi $ be a stochastic suspension of a diffeomorphism. The linear tangent
map associated to $\Phi$, and denoted by $\Phi_*$, is defined for all $X\in \Cu
(I)$ by
\begin{equation}
\Phi_* (X)=T(\Phi(X))=(\Phi (X), {\cal D} (\Phi (X))) .
\end{equation}
\end{defi}

All the quantities are well defined as diffeomorphisms send $\Cu (I)$ on $\Cu
(I)$.

\section{Invariance}

We then obtain the following notion of {\it invariance} under a one parameter
group of diffeomorphisms.

\begin{defi}
Let $\Phi=\{ \phi_s \}_{s\in \rR}$ be a one-parameter group of diffeomorphisms
and let $L$ be a functional $L: \Cu (I)\rightarrow \Cu_{\C} (I)$. The functional
$L$ is invariant under the one-parameter group of diffeomorphisms $\Phi$ if
$$L( \phi_* X )=L(X) ,\ \mbox{for all}\ \phi \in \Phi.$$
\end{defi}

As a consequence, if $L$ is invariant under $\Phi$, we have
$$L (\phi_s (X), {\cal D} (\phi_s (X)) )=L(X, {\cal D} X ) ,$$
for all $s\in \rR$ and $X\in \Cu (I)$.

\begin{rema}
We note that this notion of invariance under a one parameter group of
diffeomorphisms does not coincide with the same notion as defined by K. Yasue
(\cite{ya1}, p.332, formula (3.1)) which in our notation is given by:
$$L(\phi_s (X),\phi_s ({\cal D} X))=L(X,{\cal D} X ),\ \mbox{for all}\ s\in\rR\ \mbox{and}\ X\in \Cu (I).$$
In fact, K. Yasue definition of invariance does not reduce to the classical
notion (see for example \cite{ar},p.88) for differentiable deterministic
stochastic processes.\\

Moreover, Yasue's definition is not coherent with the invariance notion used in
his proof of the stochastic Noether's theorem (\cite{ya1},theorem 4,p.332). See
the comment below.
\end{rema}

\section{The stochastic Noether's theorem}

Noether's theorem has already been generalized a great number of times and covers
sometimes different statements \cite{kos}. Here, we follow V.I. Arnold's
(\cite{ar},p.88) presentation of the Noether theorem for
Lagrangian systems. We correct a previous work of K. Yasue
(\cite{ya1}, Theorem 4,p. 332-333).

\begin{thm}
\label{noether} Let $J_{a,b}$ be a functional on $\Cu (I)$ given by
$$J_{a,b} (X)=E\left [ \di\int_a^b L(X(t),{\cal D} X(t) ) dt \right ] .$$
with $L$ invariant under the one-parameter group $\Phi=\{ \phi_s \}_{s\in\rR}$.

Let $X\in \Cu (I)$ be a $\Cu (I)$-stationary point of $J_{a,b}$ with fixed end
points condition
$$X(a)=X_a,\ \mbox{and}\ X(b)=X_b .$$
Then, we have
$$
\di {d\over dt} E\left [ \grad_v L \left . {\partial Y \over \partial s} \right
|_{s=0} \right ] =0,$$ where
\begin{equation}
Y_s =\Phi_s (X) .
\end{equation}
\end{thm}

\begin{proof}
Let $Y(s,t)=\phi_s X(t)$ for $s\in \rR$ and $a\leq t\leq b$.

As $L$ is invariant under $\Phi=\{ \phi_s \}_{s\in \R}$, we have
$$\di {\partial \over \partial s} L(Y(s,t),{\cal D}_{\mu} Y(s,t) )=0\hskip 5mm (a.s.).$$
As $Y(.,t)$ and ${\cal D}_{\mu} Y(.,t)$ belong to $C^1 (\R )$ for all $t\in
[a,b]$ by definition \ref{group1para}, iii), we obtain
\begin{equation}
\grad_x L \cdot \di {\partial Y \over \partial s} +\grad_v L \di {\partial {\cal
D}_{\mu} Y \over \partial s} =0 \hskip 5mm (a.s.).
\end{equation}
Using (Lemma \ref{tecnicnoether},ii), this equation is equivalent to
\begin{equation}
\grad_x L \cdot \di {\partial Y \over \partial s} +\grad_v L \di {\cal D}_{\mu}
\left ( {\partial Y \over \partial s}\right ) =0 \hskip 5mm (a.s.).
\end{equation}
As $X=Y\mid_{s=0}$ is a stationary process for $J_{a,b}$, we have
\begin{equation}
\grad_x L ={\cal D}_{-\mu} \grad_v L.
\end{equation}
As a consequence, we deduce that
$$\left . \left (\left [\di {\cal D}_{\mu} \grad_v L \right ]\cdot \di {\partial Y \over \partial s} +\grad_v L \di {\cal D}_{\mu}
\left ( {\partial Y \over \partial s}\right ) \right ) \right |_{s=0}  =0
\hskip 5mm (a.s.).$$
Taking the absolute expectation, we obtain
\begin{equation}
E\left [ \left . \left (\left [\di {\cal D}_{\mu} \grad_v L \right ]\cdot \di
{\partial Y \over \partial s} +\grad_v L \di {\cal D}_{\mu} \left ( {\partial Y
\over \partial s}\right ) \right ) \right |_{s=0}  \right ] =0.
\end{equation}
Using the product rule, we obtain
$$
\di {d\over dt} E\left [ \grad_v L \left . {\partial Y \over \partial s} \right
|_{s=0} \right ] =0 ,$$ which concludes the proof.
\end{proof}

\section{Stochastic first integrals}

The previous theorem leads us to the introduction of the notion of {\it first integral} for
stochastic Lagrangian systems\footnote{Of course, one can extend this
definition to general stochastic dynamical systems.}.

\subsection{Reminder about first integrals}

Let $X$ be a $C^k$ vector field or $\R^n$, $k\geq 1$ ($k$ could be $\infty$
or $\omega$, i.e. analytic). We denote by $\phi_x (t)$ the solution of the
associated differential equation, such that $\phi_x (0)=x$ and by $S$ the set
of all
these solutions. \\

A first integral of $X$ is a real valued function $f:\R^n \rightarrow \R$ such
that for all $\phi_x (t) \in S$, we have
\begin{equation}
f(\phi_x (t) )=c_x ,
\end{equation}
where $c_x$ is a constant.\\

We have not imposed any kind of regularity on the function $f$, so that $f$ can
be just $C^0$. In this case, the existence of
a first integral does not impose many constraint on the dynamics.\\

If $f$ is at least $C^1$, then we can characterize first integrals by the
following constraint:
\begin{equation}
X\cdot f =0 .
\end{equation}

\subsection{Stochastic first integrals}

The previous paragraph leads us to searching for an analogue of the classical
notion of first integrals as a functional defined on the set of solutions of a
given stochastic Euler-Lagrange equation\footnote{Of course, this definition
will extend to arbitrary stochastic dynamical systems.} and real valued.
Looking for the stochastic Noether theorem, we choose the following definition:

\begin{defi}
Let $L$ be an admissible Lagrangian system. A functional $I:\Cu (I) \rightarrow
\R$ is a first integral for the Euler-Lagrange equation associated to $L$ if
\begin{equation}
\di {d\over dt} \left [ I(X) \right ] =0,
\end{equation}
for all $X$ satisfying the Euler-Lagrange equation.
\end{defi}

We can now interpret the stochastic Noether theorem in term of first integrals,
i.e. the fact that the invariance of the Lagrangian $L$ under of a one
parameter group of diffeomorphisms $\Phi =(\phi_s )_{s\in \R}$ induces the
existence of a first integral for the associated Euler-Lagrange equation,
defined by
\begin{equation}
I(X)=E\left [ \grad_v L \left . {\partial \phi_s X(t) \over \partial s} \right
|_{s=0} \right ] .
\end{equation}

\section{Examples}

\subsection{Translations}

We follow the first example given by V.I. Arnold (\cite{ar},p.89) for Noether theorem. Let $L$ be the Lagrangian
defined by
\begin{equation}
L(X,V)=\frac{V^2}{2}-U(X),\ \mbox{\rm where}\ X\in\R^3,
\end{equation}
$V=(V_1,V_2,V_3)\in\C^3$, $V^2:=V_1^2+V_2^2+V_3^2$ and $U$ is taken to be
invariant under the one parameter group of translations:
\begin{equation}
\phi_s(x)=x+se_1 ,
\end{equation}
where $\{e_1,e_2,e_3\}$ is the canonical basis of $\R^3$.\\

Then, by the Stochastic Noether's theorem, the quantity
\begin{equation}
E[\cal{D} X_1]
\end{equation}
is a first integral since $\partial_V L=V$ and $\partial_s\phi_s(X_1(\omega))=e_1$.

\subsection{Rotations}

We keep the notations of the previous paragraph. We consider the Lagrangian of the two-body problem in $\R^3$, i.e.
\begin{equation}
L(X,V)=q(V)-\frac{1}{|X|} \  \ \mbox{\rm  where}\ \ q(V)=\frac{V^2}{2} ,
\end{equation}
where $\mid .\mid$ denotes the classical norm on $\R^3$ defined for all $X\in \R^3$, $X=(X_1, X_2 ,X_3 )$ by
$\mid X\mid^2 =X_1^2 +X_2^2 +X_3^2$.\\

We already know that the classical Lagrangian $L$ is invariant under rotations when $X\in \R^3$ and $V\in \R^3$.
Here, we must prove that the same is true for the extended object, i.e. for $L$ defined over $\R^3\setminus \{ 0\}
\times \C^3$. This extension, as long as it is defined, is canonical. Indeed, we define $q(z)$ for $z\in \C^3$ as
\begin{equation}
q(z)=\di {1\over 2} (z_1^2 +z_2^2 +z_3^2 ) ,\ \ z=(z_1 ,z_2 ,z_3 )\in \C^3 .
\end{equation}
Note that our problem is not to discuss an analytic extension of the real valued kinetic energy but only to look for the same
function on $\C^3$ simply replacing real variables by complex one. As long as the new object is well defined this procedure
is canonical, which is not the case if we search for an analytic extension of $q$ over $\C^3$ which reduces to $q$ on
$\R^3$.\\

Our main result is then that this group of symmetry is preserved under stochastization, which is in fact a general
phenomenon that will be discuss elsewhere.

\begin{lem}
The lagrangian $L$ defined over $\R^3 \setminus \{ 0\}\times \C^3$ is invariant under rotations $\phi_{\theta,k}$ around
the $e_k$ axis by the angle $\theta$, $k=1,2,3$.
\end{lem}

The proof is based on the two following facts:
\begin{itemize}
\item As $\phi_{\theta ,k}$ is a linear map whose matrix coefficients do not depend on $t$, we have
\begin{equation}
{\cal D}_{\mu} \left [ \phi_{\theta ,k} (X) \right ] =\phi_{\theta ,k} \left [ {\cal D}_{\mu} X \right ] ,
\end{equation}
where $\phi_{\theta,k}$ is trivially extended to $\C^3$.\\

\item A simple calculation gives
\begin{equation}
\forall\ z\in\C^3,\  \ q(\phi_{\theta,k}(z))=q(z) .
\end{equation}
\end{itemize}

We easily deduce the $\phi_{\theta ,k}$ invariance of $L$, i.e. that
\begin{equation}
L(\phi_{\theta,k}X,\cal{D}(\phi_{\theta,k}X))=L(X,\cal{D}X) .
\end{equation}

We now compute: $\partial_\theta\phi_{\theta,k}(X)|_{\theta=0}=e_k\wedge X$ and
$$\partial_VL(X,\cal{D}X)\cdot \partial_\theta\phi_{\theta,k}(X)|_{\theta=0}=(X\wedge\cal{D}X)_k.$$

Therefore the expectation of the "complex angular momentum" $X\wedge\cal{D}X$
is a conserved vector ($\wedge$ is extended in a natural way to complex
vectors).

\section{About first integrals and chaotic systems}

In this section, we discuss some consequences of the stochastic Noether's
theorem in the context of chaotic
dynamical systems. The study of deterministic chaotic dynamical systems is difficult.\\

Here again, we return to the classical $n$-body problem, $n\geq 3$. In this
case, in particular for large $n$, the dynamics of the system is very
complicated and only numerical results give a global picture of the phase
space. Despite the existence of a chaotic behaviour, there exist several well
known first integrals of the system.

These integrals are used as constraints on the dynamics and can give
interesting results, as for example J. Laskar's \cite{la2} approach to the
Titus-Bode law for the repartition of the planets in the solar systems and
extra-solar systems.\\

Using our approach, we can go further by claiming that such kind of integrals
continue to exist even if one consider a more general class of perturbations
including stochasticity. We note that this result is fundamental as long as one
wants to relate numerical computations on the $n$-body problem with the real
dynamical behaviour of the solar systems, and in this particular example, the
dynamics of the protoplanetary nebulae.

\chapter{Natural Lagrangian systems and the Schr\"odinger equation}
\label{chapschrodinger}

In this section, we explore in details the stochastization procedure for
natural Lagrangian systems. In particular, by introducing a suitable analogue
of the action functional, we prove that the stochastic Euler-Lagrange equation
leads to a non-linear Schr\"odinger equation, depending on a free parameter
related to a normalization constraint. For a suitable choice of this parameter
we then obtain the classical linear Schr\"odinger equation.

\section{Natural Lagrangian systems}

In (\cite{ar},p.84), V.I. Arnold introduces the following notion of {\it
natural Lagrangian systems}:

\begin{defi}
A Lagrangian system is called natural if the Lagrangian function is equal to
the difference between kinetic and potential energy:
$$L(x,v)=T(v)-U(x) .$$
\end{defi}

As an example, we have the natural Lagrangian function associated to {\it
Newtonian mechanics}:
$$L(x,v)=\frac{1}{2}v^2-U(x),$$
where $U$ is of class $C^\infty$.

\section{Schr\"odinger equations}

\subsection{Some notations and a reminder of the Nelson wave function}

We recall that $\Lambda_d$ is the space of "good" diffusion processes. Let
$\Lambda_d^g$ be the subspace of $\Lambda_d$ whose elements have a smooth
gradient drift. We then set:
$$\cal{S}=\{X\in\Lambda_d \, \mid\, \cal{D}^2X(t)=-\nabla
U(X(t))\}.$$ For a diffusion $X$ in $\Lambda_d$ with drift $b$ and density
function $p_t(x)$, we set:
\begin{equation}
\Theta=(\R^+\times\R^d)\setminus \{(t,x),\, \mid\,  p_t(x)=0\}.
\end{equation}
If $X\in \Lambda_d^g$ then there exist real valued functions $R$ and $S$ smooth on $\Theta$ such that
\begin{equation}
\cal{D}X(t)=\left(b-\frac{\sigma^2}{2}\nabla
\log(p_t)+i\frac{\sigma^2}{2}\nabla \log(p_t)\right)(X(t))=(\nabla S+i\nabla
R)(X(t)) ,
\end{equation}
since $b$ is a gradient. Obviously:
\begin{equation}
R(t,x)=\frac{\sigma^2}{2}\log(p_t(x)).
\end{equation}
In this case, we introduce the function:
\begin{equation}
\label{wave function}
 \Psi(t,x)=e^{\di \frac{\di (R+iS)(t,x)}{K}}
\end{equation}
(where $K$ is a positive constant) called the wave function.

The wave function has the same form than that of Nelson one (see \cite{ne1}). We then set $A=S-iR$. So $\d
\Psi=e^{\frac{iA}{K}}$ and $\nabla A(t,X(t))=\overline{\cal{D}}X(t)$. For a
suitable $K$, Nelson shows that if $X$ satisfies its stochastized Newton
equation (which is the real part of ours) then $\Psi$ satisfies a Schr\"odinger
equation. We show, by using our operator $\cal{D}$, the same kind of result in
the next section.

\subsection{Schr\"odinger equations as necessary conditions}

\begin{thm}
\label{mainschrodinger}
If $X\in\cal{S}\cap\Lambda_d^g$, then the wave function (\ref{wave function})
satisfies the following non-linear Schr\"odinger equation on the set $\Theta$:
\begin{equation}
\label{non_linear}
iK\partial_t\Psi+\frac{K(K-\sigma^2)}{2}\frac{(\partial_x\Psi)^2}{\Psi}
+\frac{\sigma^2}{2}\Delta\Psi=U\Psi ,
\end{equation}
\end{thm}

\begin{proof}
As $U$ is a real valued function, $X\in\cal{S}$ implies
$$\overline{\mathcal{D}}^2X(t)=-\nabla U(X(t)).$$
The definition of $\Psi$ implies that on $\Theta$
$$ \nabla A=-iK\frac{\nabla\Psi}{\Psi}.$$
Since $\nabla A(t,X(t))=\overline{(\mathcal{D}X)(t)}$, we obtain
$$ iK\overline{\mathcal{D}}\frac{\partial_x\Psi}{\Psi}(t,X(t))=\nabla U(t,X(t)).$$
Therefore, considering the k-th component of the last equation and using lemma
\ref{darsien}, we deduce
$$
iK\left(\partial_t
\frac{\partial_k\Psi}{\Psi}+\overline{\mathcal{D}}X(t)\cdot\nabla
\frac{\partial_k\Psi}{\Psi} -i
\frac{\sigma^2}{2}\Delta\frac{\partial_k\Psi}{\Psi}\right)(t,X(t))=\partial_kU(X(t)).
$$

Now $\d \overline{\cal{D}}X(t)=-iK\frac{\nabla\Psi}{\Psi}(t,X(t))$. Thus, by
Schwarz lemma, we obtain
$$\overline{\cal{D}}X(t)\cdot \nabla\frac{\partial_k\Psi}{\Psi}=-iK\sum_{j=1}^d\frac{\partial_j\Psi}{\Psi}\partial_j\frac{\partial_k\Psi}{\Psi}=-\frac{iK}{2}\partial_k\sum_{j=1}^d\left(\frac{\partial_j\Psi}{\Psi}\right)^2,$$
and
$$\Delta \frac{\partial_k\Psi}{\Psi}=\sum_{j=1}^d\partial_j^2\frac{\partial_k\Psi}{\Psi}=\partial_k\sum_{j=1}^d\partial_j\left(\frac{\partial_j\Psi}{\Psi}\right)=\partial_k\sum_{j=1}^d\frac{\partial_j^2\Psi}{\Psi}-\left(\frac{\partial_j\Psi}{\Psi}\right)^2.$$

Therefore
$$iK\partial_k\left(\frac{\partial_t\Psi}{\Psi}+i\frac{\sigma^2-K}{2}\partial_k\sum_{j=1}^d \left(\frac{\partial_j\Psi}{\Psi}\right)^2
-i
\frac{\sigma^2}{2}\frac{\Delta\Psi}{\Psi}\right)(t,X(t))=\partial_kU(X(t)).$$
By adding an appropriate function of $t$ in $S$, we can arrange the constant in
$x$ of integration in equation to be zero, and formula (\ref{non_linear})
follows as claimed.
\end{proof}

In order to recover the classical linear Schr\"odinger equation, we must choose
the normalization constant $K$. The main point is that in this case, we obtain
a clear relation between the modulus of the wave function and the density of
the underlying diffusion process. Precisely, we have:

\begin{cor}
We keep the notations and assumptions of theorem (\ref{mainschrodinger}). We
assume that
$$K=\sigma^2.$$
Then the wave functional $\Psi$ satisfies the linear Schr\"odinger equation
\begin{eqnarray}
i\sigma^2\partial_t\Psi+\frac{\sigma^4}{2}\Delta\Psi=U\Psi ,
\end{eqnarray}
Moreover, if $p_t(x)$ is the density of the process $X(t)$ at point $x$, then
we have
$$(\Psi\overline{\Psi})(t,x)=p_t(x) .$$
\end{cor}

\begin{proof}
$K=\sigma^2$ kills the non-linearity in equation (\ref{non_linear}) and
furthermore
$$ \log(\Psi\overline{\Psi})=\frac{2}{K}R =\frac{2}{\sigma^2}R=\log(p).$$
which concludes the proof.
\end{proof}

\subsection{Remarks and questions}

\begin{itemize}
\item Obviously $\Lambda_1\subset\Lambda_1^g$ since $b$ is continuous.\\

\item A natural question is to know if the converse of the corollary of (\label{mainschrodinger}) is true. More
precisely, if $\Psi$ satisfies a linear Schr\"odinger equation, can we construct a process $X$ which belongs to
$\cal{S}\cap\Lambda_d^g$ and whose density is such that $p_t(x)=|\Psi(t,x)|^2$ ?\\

R. Carmona tackled the problematic of the so-called Nelson processes and proved in \cite{car} under some conditions
the existence of a process $X$ with gradient drift related to $\Psi$ and whose density is such that $p_t(x)=|\Psi(t,x)|^2$.
However we do not know if this process belongs to our space of good diffusions processes (which may turn to be a little
restrictive class in this case), but we can prove formally, \textit{i.e.} even so assuming that the formulae of the
stochastized derivative to a function of the process holds, that $X$ satisfies the Newton stochastized equation.
Therefore, this leads one to question the extension of the derivative operator and the way it acts on a large class
of processes. This problem will be treated in a forthcoming paper (See \cite{ccbd1}).\\

\item The fact that a process $X$ satisfies the stochastized Newton equation of Nelson implies $(D^2-D_*^2)X=0$ (for the
potential $U$ is real). This is a general fact for diffusion with gradient drift. Indeed, we can prove:

\begin{lem}
Let $X\in\Lambda_d$, $b$ its drift and $p$ its density function. Let $G_i$ be the i-th column of the
matrix $(G_{ij}):=(\partial_jb_i-\partial_ib_j)$. Then $(D^2-D_*^2)X=0$ if and only if for all $t>0$, ${\rm div}(p_t G_i)=0$.
\end{lem}

Thus, if $X\in\Lambda_d^g$ it is clear that $(D^2-D_*^2)X=0$ since the form $\sum b_k\partial_k$ is closed and so $G=0$.
An interesting question is then to know if the converse is true. So we may wonder ourselves if
$\cal{S}\subset\Lambda_d^g$.\\

The difficulty relies on the fact that $p$ and $b$ are related via the Fokker-Planck equation, so the condition
${\rm div} (p_t G_i)=0$ may not be the good formulation. However, one could use the work of S. Roelly and M. Thieullen
in \cite{roelly_thieullen} who use an integration by parts via Malliavin Calculus to characterize gradient diffusion,
in order to give a positive or negative answer to our question.\\

\item A basic notion in mechanics is that of {\it action} (see
\cite{ar},p.60). The action associated to a Lagrangian system is in general
obtained via the action functional. In our framework, a natural definition for
such an action functional is given by:

\begin{defi}
Let $\mathcal{A}$ be the functional defined on $[a,b]\times
\EuScript{C}^1([a,b])$ by:
\begin{equation}
\forall t\in[a,b],\ \forall X\in\EuScript{C}^1(I),\
\mathcal{A}(t,X)=E\left[\int_a^t
L(X_s,(\mathcal{D}X)_s)ds\left|X_t\right.\right].
\end{equation}
This functional is called the {\it action functional}.
\end{defi}

Using this action functional, we have some freedom to define the corresponding ``action". The natural
one is defined by
\begin{equation}
A_X(t,x)=E\left[\int_a^t L(X_s,(\mathcal{D}X)_s)ds\left|X_t=x\right.\right] .
\end{equation}
Usually, the wave function associated to $A_X$ an denoted by $\tilde{\psi}$ is then defined as
\begin{equation}
\tilde{\psi}_X (t,x)=\exp^{i A_X (t,x)} .
\end{equation}
However, it is not at all clear that such kind of function satisfies the {\it gradient condition}, i.e. that
\begin{equation}
\label{gradientaction}
\nabla A (t,X(t)) =\overline{\cal D} X (t) ,
\end{equation}
which is fundamental in our derivation of the Schr\"odinger equation.\\

However, the condition \ref{gradientaction} is equivalent to prove that the real part of ${\cal D} X$ is a gradient,
which is not at all trivial in dimension greater than two.
\end{itemize}

\section{About quantum mechanics}

Even if we look for dynamical systems, our work can be used in the context of
the so-called {\it Stochastic mechanics}, developed by Nelson \cite{ne1}. The
basic idea is to reexpress quantum mechanics in terms of random trajectories.
We refer to \cite{cz} for a review.\\

The stochastic embedding theory can be seen as a quantization procedure, i.e. a
formal way to go from classical to quantum mechanics. This approach is already
different from Nelson's approach, which do not define a rigid procedure to
associate to a given equation a stochastic analogue. Moreover, the acceleration
defined by Nelson as
\begin{equation}
a(X)=\di {DD_* (X)+D_* D(X)\over 2} ,
\end{equation}
is only a particular choice. Many authors have tried to justify this form
(\cite{pav1},\cite{pav2}) or to try another one. In our context, the form of
the acceleration is fixed and corresponds, as in the usual case, to the second
(stochastic) derivative of $X$. As a consequence, stochastic embeddings can be
used to provide a conceptual framework to stochastic mechanics. We refer to
\cite{pav1} where a complex valued velocity for a stochastic process is
introduced
corresponding to the stochastic derivative of $X$.\\

However, stochastic mechanics as well as its variants have many drawbacks with
respect to the initial wish to describe quantum mechanical behaviours. We refer
to \cite{ne3} and \cite{cz} for details. This is the reason why we will not
develop further this topic.

\chapter{Stochastic Hamiltonian systems}
\label{hamilton}
\setcounter{section}{0}
\setcounter{equation}{0}

In this part, we introduce the stochastic pendant of Hamiltonian systems for classical Lagrangian systems.
The strategy is first to define the stochastic analogue of the classical momentum. We then define a
stochastic Hamiltonian. However, this Hamiltonian is not obtained by the classical stochastic embedding
procedure. This is due to the fact that the momentum process is complex valued. As a consequence, we must
modify the procedure in order to obtain a coherent picture between the classical formalism and the stochastic
one. This leads us to define the stochastic Hamiltonian embedding procedure which reflects in fact the non
trivial character of the underlying stochastic symplectic geometry to develop. Having the stochastic Hamiltonian
we prove a Hamilton least action principle using our stochastic calculus of variations. We then obtain an analogue of
the Lagrangian coherence lemma in this case up to the fact that the underlying stochastic embedding procedure is now
the Hamiltonian one.

\section{Reminder about Hamiltonian systems}

We denote by $I$ an open interval $(a,b)$, $a<b$.\\

Let $L:\R^d \times \R^d \times \R \rightarrow \R$ be a convex Lagrangian. The
Lagrangian functional over $C^1 (\R )$ is defined by
\begin{equation}
L: \ \left .
\begin{array}{lll}
C^1 (\R ) & \longrightarrow & C^1 (\R ), \\
x & \longmapsto & L(x,\dot{x} , t) .
\end{array}
\right .
\end{equation}

We can associate to $L$ a Hamiltonian function using the {\it Legendre transformation} (\cite{ar},p.65). From the
functional side, this induces a change of point of view, as the functional is not seen as acting on $x(t)$,
which is the so-called {\it configuration space} of classical mechanics, but on
$(x(t),\dot{x} (t))$ which is associated to the {\it phase-space}. This dichotomy between position and velocities has of course many
consequences, one of them being that the system is more symmetric (the
symplectic structure).

\begin{defi}
Let $L(x,v)$ be an admissible Lagrangian system. For all $x\in C^1$, we denote by
\begin{equation}
p (x)=\di {\partial L\over \partial v} (x,\dot{x} ),
\end{equation}
the momentum variable.
\end{defi}

We now introduce an important class of Lagrangian systems.

\begin{defi}
Let $L(x,v)$ be an admissible lagrangian system. The Lagrangian $L$ is said to possess the Legendre property if
there exists a function $f:\R^d \rightarrow \R^d$, called the Legendre transform, such that
\begin{equation}
\dot{x}=f(x,p) ,
\end{equation}
for all $x\in C^1$.
\end{defi}

Most classical examples in mechanics possess the Legendre property. This follows from the convexity of $L$ in the second variable
(see \cite{ar},p.61-62).\\

We can introduce the fundamental object of this section:

\begin{defi}
Let $L$ be an admissible Lagrangian system which possesses the Legendre property. The Hamiltonian function associated to
$L$ is defined by
\begin{equation}
H(p,x)=pf(x,p) -L(x,f(x,p)) ,
\end{equation}
where $f$ is the Legendre transform.
\end{defi}

The Hamiltonian function plays a fundamental role in classical mechanics. We introduce the stochastic analogue in the next
section.

\section{The momentum process}

A natural stochastic analogue of the momentum variable is defined as follow:

\begin{defi}
\label{stocmomentum}
Let $L(x,v)$ be an admissible Lagrangian system. For all $X\in \Cu (I)$, we
define the stochastic process $P(t)$, called the canonical momentum process, by
\begin{equation}
\label{defimomentum}
P (t)=\di {\partial L\over \partial v} (X(t),{\cal D}X (t)).
\end{equation}
\end{defi}

This definition can be made more natural using the embedding $\iota$ defined from $C^0 (I)$ on ${\cal P}_{\rm det}$ and
the linear tangent map introduced in chapter \ref{noetherchap}. Indeed, the momentum process can be viewed as a functional on
$X\in \Cu (I)$, $P:\Cu (I) \rightarrow {\cal P}_{\C}$ defined by (\ref{defimomentum}). We have for all
$X\in {\cal P}_{\rm det}^1 =\iota (C^1 (I))$,
\begin{equation}
P (X) =\iota (p (x)),
\end{equation}
where $x\in C^1 (I)$ is such that $X=\iota(x)$. As by definition, we have
\begin{equation}
\label{formmomen}
\iota (p (x))=p (\iota (x))=p(X) .
\end{equation}
As $p$ keeps a sense for $X\in \Cu (I)$, we extend formula (\ref{formmomen}) to $\Cu (I)$ leading to definition
\ref{stocmomentum}.\\

If we assume that the Lagrangian possesses the Legendre property, then there exists a Legendre transform $f$ such that
for all $x\in C^1$, $\dot{x}=f(x,p)$. We can ask if such a property is conserved for the momentum process. We have:

\begin{lem}
Let $L(x,v)$ be an admissible Lagrangian system possessing the Legendre property. Let $f$ be the Legendre transform
associated to $L$. We have
\begin{equation}
{\cal D} X(t)=f(X,P) ,
\end{equation}
for all $X\in \Cu (I)$.
\end{lem}

We can now define the stochastic Hamiltonian associated to $L$:

\begin{defi}
\label{stochamilton}
Let $L(x,v)$ be an admissible Lagrangian system possessing the Legendre property. The stochastic
Hamiltonian system associated to $L$ is defined by
\begin{equation}
H:\ \left .
\begin{array}{lll}
{\cal P}_{\C} \times \Cu (I) & \longrightarrow & {\cal P}_{\C} \\
(P,X) & \longmapsto & P f(X,P)-L(X,f(X,P) ) .
\end{array}
\right .
\end{equation}
\end{defi}

\section{The Hamiltonian stochastic embedding}

As in the previous chapter, we want to use the stochastic embedding procedure to associate a natural
stochastic analogue of the Hamiltonian equations. However, we must be careful with such a procedure, as already
discussed in chapter \ref{chapembed}, $\S$.\ref{secondorder}. Indeed, the embedding procedure does not allow us to fix the notion of
embedding for {\it systems} of differential equations. Moreover, we must keep in mind that the principal idea behind the
Hamiltonian formalism is to work not in the {\it configuration space}, i.e. the space of positions, but in the {\it phase
space}, i.e. the space of positions and momenta. As the stochastic speed is by definition complex, this induces a
particular choice for the embedding procedure in the case of Hamiltonian differential equations.

\begin{defi}
Let $F:\R^d \times \C^d \mapsto \C$ be a holomorphic function, real valued on real arguments. This function defines a
real valued functional over $C^1 (I)\times C^1 (I)$, for $I$ a given open interval of $\R$. The
Hamiltonian embedding of the functional $F$ is the functional denoted by $F_S$, defined on $\Cu (I) \times
{\cal P}_{\C} (I)$ by $H$, i.e.
\begin{equation}
F_S (X,P) (t) =F(X (t),P (t)) .
\end{equation}
\end{defi}

We denote by $S_H$ the procedure associating the stochastic functional $F_S$ to $F$. This procedure reduces to
change the functional spaces for $F$ from $C^1 (I) \times C^1 (I)$ to $\Cu (I) \times {\cal P}_{\C}$.\\

The main property of the Hamiltonian stochastic embedding procedure (and in fact it can be used as a definition) is to lead
to a coherent definition with respect to the momentum process. Precisely, we have:

\begin{lem}[Legendre coherence lemma]
Let $L(x,v)$ be an admissible Lagrangian system possessing the Legendre property. The following diagram commutes
\begin{eqnarray}
\xymatrix{
  & (x,p) \ar[d]_{S_H} \ar[r]^{H} & H(x,p) \ar[d]^{S_H}       \\
  & (X,P) \ar[r]_{H_S}   & H (X,P)
  }
\end{eqnarray}
\end{lem}

The proof follows essentially from the fact that the stochastic Hamiltonian embedding of the functional $H$, denoted by
$H_S$ coincide with the definition \ref{stochamilton} of the stochastic Hamiltonian system associated to $H$ via the Legendre transform and
the definition of the momentum process.

\section{The Hamiltonian least action principle}

Using the stochastic Hamiltonian function, we can use the stochastic calculus of variations in order to obtain the
set of equations which characterize the stationary processes of the following functional:
\begin{equation}
I_{a,b} (X,P)=E \left [
\int_a^b (P(t){\cal D} X -H(X(t),P(t)))\, dt \right ] ,
\end{equation}
defined on $\Cu (I)\times {\cal P}_{\C}$.\\

In order to apply our stochastic calculus of variations, we restrict our attention to $I$ on
$\Cu (I)\times \Cu (I)$. The fundamental result of this section is the following:

\begin{thm}
A necessary and sufficient condition for an $L$-adapted process $(X,P)$ to be ${\cal N}^1 (I)$-stationary process of
the functional $I_{a,b}$ with fixed end points $(X(a),P(a))=(X_a ,P_a ) \in H$, $(X(b),P(b))=(X_b ,P_b) \in H$ is that
it satisfies the stochastic Hamiltonian equations
\begin{equation}
\left .
\begin{array}{lll}
{\cal D} X & = & \di {\partial H\over \partial P} (X(t),P(t)),\\
{\cal D} P & = & -\di{\partial H\over \partial X} (X(t),P(t)) .
\end{array}
\right .
\end{equation}
\end{thm}

\begin{proof}
We must use the weak least action principle using the process $Z=(X,P)\in \Cu (I )\times {\cal P}_{\C}$ and the
Lagrangian denoted by ${\cal L}$ defined on $\R^d \times \C^d \times \C^d \times \C^d$ by
\begin{equation}
{\cal L} (x,p ,v ,w)=pv  -H(x,p) .
\end{equation}
As ${\cal L} (x,p,v,w)=L(x,v)$ formally via the Legendre transform, and $L$ is assumed to be admissible, we deduce that
${\cal L}$ is again admissible.\\

Let $\delta Z$ be a ${\cal N}^1 (I)$ variation of the form $Z+\delta Z=(X+X_1 ,P+P_1 )$, where $X_1$ and $P_1$ are ${\cal N}^1$
processes.

The Euler-Lagrange equation associated to ${\cal L}$ is given by
\begin{equation}
\left .
\begin{array}{lll}
\di {\partial {\cal L} \over \partial x} (Z(t), {\cal D} Z (t)) - {\cal D}_{\mu} \left [
\di {\partial {\cal L} \over \partial v} (Z(t), {\cal D} Z (t)) \right ] =0,\\
\di {\partial {\cal L} \over \partial p} (Z(t), {\cal D} Z (t)) - {\cal D}_{\mu} \left [
\di {\partial {\cal L} \over \partial w} (Z(t), {\cal D} Z (t)) \right ] =0.
\end{array}
\right .
\end{equation}
An easy computation leads to
\begin{equation}
\left .
\begin{array}{lll}
-\di {\partial H\over \partial x} (Z(t),{\cal D} Z(t)) -{\cal D}_{\mu} P (t) & = & 0 ,\\
{\cal D} X (t) - \di {\partial H\over \partial p} (Z(t), {\cal D} Z(t)) & = & 0 .
\end{array}
\right .
\end{equation}
This concludes the proof.
\end{proof}

\begin{rema}
In this proof we do not need a uniform assumption on the set of variations as the Lagrangian does not depend on the
variable $w$. In fact, we can assume a variation in the direction $P$ which belongs to $\Cu (I)$.
\end{rema}

\section{The Hamiltonian coherence lemma}

In this section, we derive the Hamiltonian analogue of the Lagrangian coherence lemma.

\begin{lem}[The Hamiltonian cohrence lemma]
Let $H: \rR^d \times \rR^d \rightarrow \rR$ be an admissible Hamiltonian system. Then, the following diagram commutes
\begin{eqnarray}
\xymatrix{
  & H(x(t),p(t)) \ar[d]_{\mbox{\rm Least action principle}} \ar[r]^{S_H} & H(X(t),P(t)) \ar[d]^{
  \mbox{\rm Stochastic least action principle}}       \\
  & (HE) \ar[r]_{S_H}   & (SHE)               }
\end{eqnarray}
\end{lem}

The main point is that this result is not valid if one replaces the Hamiltonian stochastic embedding by the natural stochastic
embedding that we have used up to now. We can keep the classical embedding procedure only when dealing with real valued versions
of the stochastic derivative. For example, if one deals with the reversible stochastic embedding procedure, we obtain a unified
stochastic embedding procedure for both Lagrangian an Hamiltonian systems. We think however that as well as the complex nature of
the stochastic derivative has a fundamental influence on the form of the stochastic Lagrangian equations, i.e. that we
obtain the Nelson acceleration, the fact to move from $S$ to $S_H$ reflects a basic properties of the underlying
stochastic symplectic geometry we must take into account this complex character of the speed. This problem will be studied in
another paper.

\chapter{Conclusion and perspectives}
\label{conclusion}

This part aims at discussing possible developments and applications of
the stochastic embedding procedure.

\section{Mathematical developments}

\subsection{Stochastic symplectic geometry}

The Hamiltonian formalism developed in the last part suggest the introduction
of what can be called a stochastic symplectic geometry. An interesting
construction of symplectic structures on Hilbert spaces is given in \cite{kuk}.\\

The main point here is to construct an analogue of the geometrical structure which puts in
evidence the very particular symmetries of the Lagrangian equations in classical mechanics. There
exists already many attempt to construct a given notion of symplectic geometry or at least a given geometry for stochastic
processes, but they are as far as we know of a different nature. We refer to the book of Elworthy, LeJan and Li \cite{ell}
for an overview. These geometries are only associated to stochastic processes and translate into data of geometrical
nature properties of the underlying stochastic processes (like the Riemannian or sub-Riemannian structure associated to
Brownian motions and diffusions).\\

A recent work of J-C. Zambrini and P. Lescot (\cite{zl1} and \cite{zl2}) deals specifically with symplectic geometry and
a notion of integrability by quadratures. \\

For a discussion of integrability in our context see section \ref{pdestoc}.

\subsection{PDE's and the stochastic embedding}
\label{pdestoc}

The stochastic embedding of Lagrangian systems over diffusion processes lead to
a PDE governing the density of the solutions of the stochastic Euler-Lagrange
equation. Moreover, we have defined a stochastic Hamiltonian system naturally
associated to the Lagrangian. However, some classical PDEs, as for example the
Schr\"odinger equation, possess an Hamiltonian formulation. This remark, which
goes back to the work of Zakharov V.E. and Faddeev D. \cite{zd} is now an
important subject in PDEs known as Hamiltonian PDEs (see for example
\cite{kuk}). As a consequence, we have the following situation:
\begin{equation}
\left .
\begin{array}{ccc}
H_S & & \\
\downarrow & & \\
PDE & \longrightarrow & H
\end{array}
\right .
\end{equation}
Of course the relation between the PDE and $H_S$ is not of the same nature as the relation with $H$.\\

In the sequel, we list a number of problems and questions which naturally arise from the previous diagram:\\

\begin{itemize}
\item There exists a notion of completely integrable Hamiltonian PDE (see \cite{kuk}). What about out stochastic
Hamiltonian systems ?
\end{itemize}

Assuming that we have a good notion of integrability for $H_S$, we have the
following questions:\\

\begin{itemize}
\item Are there any relations between the integrability of $H$ and $H_S$?
\item Is there a stochastic analogue of the Arnold-Liouville theorem?
\item Is there a special set of ``coordinates" similar to the action/angle variables? \\
\end{itemize}

We note that there already exists such a notion for Hamiltonian PDEs (see
\cite{zd}).

\begin{itemize}
\item Is there a notion of integrability by ``quadratures"?
\end{itemize}

In that respect, we think about Lax work \cite{lax} on the integrability of PDEs.

\section{Applications}

\subsection{Long term behaviour of chaotic Lagrangian systems}

The dynamical behaviour of unstable or chaotic dynamical systems is far from
being understood, unless we restrict to a very particular class of systems like
hyperbolic systems or weak version of hyperbolicity. This question arises
naturally for small perturbations of Hamiltonian systems for which there exists
a large family of results dealing with this problem, as for example the KAM
(Kolmogorov-Arnold-Moser) theorem, Nekhoroshev theorem and special phenomena
like the Arnold diffusion related to
the so-called quasi-ergodic hypothesis.\\

Unfortunately, these results are difficult to use in concrete situations and
only direct numerical simulations provide
some understanding of the dynamics \cite{ladu}.\\

There exists of course ergodic theory which tries to look for weaker information
on the dynamics than a direct qualitative approach. However, this
theory leads also to very difficult problems when one tries to implement it, as
for example in the case of Sina\"i billiard. Moreover, there is a widely
opinion in the applied community that the long term behaviour of a chaotic
systems is more or less equivalent to a stochastic process. One example of such
opinion is well expressed in the article of J. Laskar \cite{la2} in
the context of the chaotic behaviour of the Solar system: ``Since the
characteristic time scale for the divergence of nearby orbits in the Solar
system is approximately 5 Myr, the orbital evolution of the planet becomes
practically unpredictable after 100 Myr. Thus in the long term, the motion of
the Solar system may be described by a random process, where orbits wander
erratically in a chaotic
zone."\\

What are the arguments leading to this idea ? \\

The first point is that chaotic dynamical systems are in general characterized
by the so-called sensitivity to initial conditions, meaning that a small error
on the initial condition leads to very different solutions. Of course, one must
quantify this kind of sentence, and we can do that, with more or less
canonicity, by introducing Lyapounov exponents and Lyapounov time. Whatever we
do, there is a non canonical data in this, which is precisely to what
extent we consider that two solutions are different. This must be a matter of
choice for a given system, and cannot be fixed by any mathematical tool. In the
sequel, we assume that a system is sensitive to initial conditions in some
region $R$ of the phase space, and for a given metric, if for all $x_0 \in R$
and all $\epsilon >0$, the distance at time $t$ between a trajectory starting
at $x_0$ and $x_0 +\epsilon$, denoted by $d(t)$ is \footnote{As we already
stress, we can in some situations gives a precise meaning to all this point,
like for example in the Smale Horseshoes, but this is far to cover the wide
variety of chaotic behaviour which are studied in the applied literature.}
approximately given by
\begin{equation}
d(t)=\epsilon e^{t/T} ,
\end{equation}
where $T>0$ is the so-called Lyapounov time or horizon of predictability for
the system\footnote{In concrete systems, one must involve a macroscopic scale
(see \cite{dou},p.17), which bound the admissible size of an error on a
prediction. Here, this quantity is arbitrary replaced by $e$.}. For an example
of such an estimate, we refer to J. Laskar \cite{la3} where he gives numerical
evidences for
the chaotic behaviour of the solar system.\\

As a consequence, for $t$ sufficiently large with respect to $T$, we have no
prediction any more, or in other words, we can not assign to a given prediction a
precise initial condition. We then have lost the deterministic character of
the equations of motions. An idea is then to say that one musts then consider
not a fixed initial condition $x_0$, but a given random variable representing
all the possible behaviours (kind of trajectories) one is lead to after a fixed
time $t$: for example, $\epsilon >0$ being fixed, we consider all the
intersections of trajectories starting in the disk $D(x_0 ,\epsilon )$ with the
ball $B(x_0 ,\epsilon )$. We then obtain a family of directions. Assuming that
we can compute an average over the family of such a quantity which obtain an
averaged direction which select a given point of the ball $B(x_0 ,\epsilon )$.
We then follow the selected trajectory during the time $t$, and continue again
this procedure. Such a construction is reminiscent of the classical
construction of the Brownian motion (see \cite{karatzas},p.66). Of course, this
programme can only be carried in some specific examples. We refer to the article
of
Y. Sina\"i \cite{sin1} for an heuristic introduction to all these problems.\\

If we agree with the previous heuristic idea, one can then ask for the
following: how is the underlying stochastic process governed by the dynamical system ?\\

We return again to the Hamiltonian/Lagrangian case. The stochastic embedding
procedure answers precisely this question. The stochastic Euler-Lagrange
equation is the track of the underlying Lagrangian system on stochastic
processes. As a consequence, we can think that we are able to capture even the
desired long term behaviour of the Lagrangian system using this
procedure.\\

In order to support our point of view, we suggest the following strategy:\\

Consider a perturbation of a completely integrable Hamiltonian system
$H_{\epsilon} (x)=h(x)+\epsilon f(x)$, with $x\in \R^{2n}$ for example. Let us
assume that $h(x)$ leads to a particular PDE under stochastic embedding, which
can be well understood and solved. The long term behaviour of the completely
integrable Hamiltonian system is trivial. This not the case for the stochastic
analogue. What about the long term behaviour of $H_{\epsilon}$ ? We think that
it is controlled by the stochastic analogue of the unperturbed Hamiltonian.
This result is related to a kind of stochastic stability which we must define.
However, this approach can be tested on a wide variety of examples, in
particular celestial mechanical problems.

\subsection{Celestial mechanics}

There exist many theories dealing with the problem of the formation of
gravitational structures. For planetary systems this question is related to a
long standing problem related to the ``regular" spacing of planets in the Solar
system. This problem which goes back to Kepler (1595), Kant (1755), von Wolf
(1726), Lambert (1761), takes a mathematical form under the Titius (1766)
formulation of the so called Titius-Bode law giving a geometric progression of
the distance of the planets from the sun. We refer to the book of Nieto
\cite{nieto} for more details. Even if this empirical law fails to predict
correctly the real distance for the Planet Pluto for example, its interest is
that it suggests that the repartition of exoplanet orbital semi-major axes
could satisfy a simple law. As a consequence, one searchs for a possible
physical/dynamical theory supporting the existence of such kind of law.
Moreover, the discovery of many exo-planetary systems can be used to test
if the theory is based on universal phenomena and not related to our knowledge of the Solar system.\\

All the actual theories about the origin of the solar system presuppose the
formation of a protoplanetary nebula, formed by some material (gas,
dust, etc ...) with a central body (a star or a big planet). We refer to
Lissauer \cite{lissau}
for more details.\\

Instead, we use a simplified model consisting of a large central body of mass
$m_0$ with a large number of small bodies $(m_j )_{j=1,\dots ,n}$, whose mass
is assumed to be small with respect to $m_0$. The main problem is to understand
the long term dynamics of this  model.\\

Following the work of Albeverio S., Blanchard Ph. and R. Hoegh-Krohn
(\cite{abh}, see also \cite{abh2}), we can modelize the motion of a given grain
in the protoplanetary nebula by a stochastic process (see
\cite{abh},p.366-367), more precisely a diffusion process. The problem is then
to find what is the equation governing the dynamics of such a stochastic
process. Using our stochastic embedding theory, we can use the classical
formulation in order to obtain the
desired equation. This question will be detailed in a forthcoming article.\\

The main idea behind stochastic modelisation is the following:\\

The motion of a given small body in a protoplanetary nebula is given by the
Kepler model and a perturbation due to the large number of number of small
bodies. In \cite{abh}, this perturbation is replaced by a white noise. As a
consequence, the movement of a small body is assumed to be described by a
diffusion process. It must be noted that this assumption is related to a number
of arguments, one of them being that the dynamics of the underlying classical
system is unstable. We then return to our previous description of the chaotic
behaviour of a dynamical system. However, using the stochastic embedding theory,
we can try to justify the passage from a classical motion to a stochastic one looking at the following problem:\\

Let $L_{\epsilon} =L_{\rm Kepler} +P_{\epsilon}$, be the Lagrangian system
describing the dynamics of our model. The Lagrangian $L_{\rm Kepler}$ is the
classical Lagrangian of the Kepler problem, and $P_{\epsilon}$ is the
perturbation. Using the stochastic embedding theory, we can deduce two
stochastic dynamical systems, one associated to $L_{\epsilon}$ and denoted by
$S_{\epsilon}$ and one associated to $L_{\rm Kepler}$ denoted by $S_{\rm
Kepler}$. If the previous strategy to replace the perturbative effect by a
White noise is valid, then we must have a kind of stochastic stability between
$S_{\rm Kepler}$ and $S_{\epsilon}$. The notion of stochastic stability must be
defined rigorously and be consistent with the stochastic embedding
theory\footnote{It must be noted that there exists already several notion of
stochastic stability in the literature, as for example Has'inskii \cite{has},
Kushner \cite{kush} and more recently Handel \cite{han}.}. Why such a stability
result is reasonable ? The main thing is that we already look in $S_{\rm
Kepler}$ for statistical properties of the set of trajectories of stochastic
(diffusion) processes under the Kepler Lagrangian. There is no reason that the
statistic of this trajectories really differs when adding a small perturbation.
This is of course different if one look for the underlying deterministic
system. All these questions will be studied in a forthcoming paper.

\subsection{Strange attractors}

{\it Strange attractors} play a fundamental role in turbulence and lead to many
difficult problems. Most of the time, one is currently interested in the
geometrical properties of attractors (Hausdorf dimension,...), special
dynamical properties (existence of an SRB (Siba\"i-Ruelle-Bowen) measure
\cite{vian}, stability under perturbations....). However, focusing on a given
attractor hides the fact that most of the time we can not predict from the
equation the {\it existence} of such an attractor. This is in particular the
case for the {\it Lorenz attractor} or the {\it Henon attractor}. These
attractors are obtained numerically. In some models, we can construct a
geometric model from which we can prove the existence of such a structure (this
is the case for the geometric Lorenz model) \cite{guck}. For example, S. Smale
\cite{sma} asks for an existence proof for the Lorenz equation of the attractor.
This has been done recently by W. Tucker (\cite{tuc2},
\cite{tuck3}). However, no general strategy exists in order to predict such an attractor.\\

Our idea is to use the stochastic embedding theory in order to predict the
existence of such an object. Let us consider the Lorenz equations. These
equations are not a Lagrangian system. However, there exits a canonical
embedding in a Lagrangian system (see the report of M. Audin \cite{aud}). This
lagrangian can then be studied via the stochastic embedding procedure. The
solutions are stochastic processes whose density is controlled by a PDE. As we
already explain, we expect that the long term behaviour of the system is coded
by this PDE. As the long term dynamics of the Lorenz system if precisely
supported by the Lorenz attractor, we think that this structure can be detected
in the PDE (as a stationary state for
example).\\

We can also take this problem as a first step towards understanding the
existence of coherent structures in chaotic dynamical systems. Moreover, the
Lorenz attractor is widely studied and there exists a great amount of results
like the existence of a unique SRB measure (see \cite{tuck3}). We can then take
this example as a good system to compare classical methods of {\it ergodic
theory} and our approach. For more problems related to the Lorenz attractor,
SRB measure $\dots$, see (\cite{vian2},\cite{vian3}).

\chapter*{Notations}

\ni $d$: dimension\\

\ni $(\Omega ,{\cal A}, P)$ a probability space
\vskip 5mm
\ni - {\bf Stochastic processes}\\

\begin{itemize}
\item We denote by
$$
dX =b(t,X) dt +\sigma (t,X) dW ,\eqno{(*)}
$$
the stochastic differential equation where $b$ is the drift, $\sigma$ the diffusion matrix and $W$ is a $d$-dimensional
Wiener process defined on $(\Omega ,{\cal A}, P)$.

\item We denote by $X(t)$ the solution of (*) and by $p_t (x)$ its density (when it exists) at point $x$.
\item $\sigma (X_s ,a\leq s\leq b )$: the $\sigma$-algebra generated by $X$ between $a$ and $b$
\item ${\cal F}_t$: an increasing $\sigma$ algebra
\item ${\cal P}_t$: an decreasing $\sigma$ algebra
\item $E\left [ \bullet \mid {\cal B} \right ]$: the conditional expectation.
\item $\parallel .\parallel$: norm on stochastic processes.
\end{itemize}

\vskip 5mm
\ni - {\bf Functional spaces}\\

\begin{itemize}
\item ${\cal P}_{\R}$: real valued stochastic processes
\item ${\cal P}_{\C}$: complex valued stochastic processes
\item ${\cal P}_{\rm det}$: the set of deterministic stochastic processes
\item ${\cal P}_{\rm det}^k$: the set of deterministic stochastic processes such that $X(\omega )$ is of class $C^k$
\item $\Lambda_d$: good diffusion processes
\item $\Lambda_d^g$: good diffusion processes with a gradient drift
\item $L^p (\Omega )$: set of random variables which belongs to $L^p$
\item ${\bf L}^2$: the set of real valued processes which are ${\cal P}_t$ and ${\cal F}_t$ adapted and
such that $E\left [ \di\int_0^1 X_t^2 \, dt \right ] <\infty$.
\item $C^{1,2} ((0,1)\times \R^d )$ the set of function which are $C^1$ in the first variable and $C^2$ in the second
one.
\item ${\cal N}^1$: the set of Nelson differentiable processes.
\end{itemize}
\vskip 5mm
- {\bf Operators}\\

\begin{itemize}
\item $\nabla$: the gradient
\item $\Delta$: the Laplacian
\item Let $f(x_1 ,\dots ,x_n )$ be a given function. We denote by $\partial_{x_i} f$ the partial derivative of $f$ with respect
to $x_i$
\item Let $f(x_1 ,\dots ,x_n ,y_1 ,\dots ,y_m )$ be a given function. We denote by $\partial_x f$, $x=(x_1,\dots ,x_n )$ the
partial differential of $f$ in the direction $x$.
\item $D$: Nelson forward derivative
\item $D_*$: Nelson backward derivative
\item $\cal D$: the stochastic derivative
\item $D^n$, $D_*^n$, ${\cal D}^n$: the $n$-th iterate of $D$, $D_*$ or ${\cal D}$
\item $\bab$ and $\bab_*$: adapted forward and backward derivative
\end{itemize}
\vskip 5mm
$k\geq 1$
\vskip 3mm
\begin{itemize}
\item ${\cal C}^k$: the set of real valued processes which are ${\cal P}_t$ and ${\cal F}_t$ adapted and
such that ${\cal D}^i$ exists, $1\leq i\leq k$.
\item ${\cal C}^k_{\C}$: the set of complex valued processes which are ${\cal P}_t$ and ${\cal F}_t$ adapted and
such that ${\cal D}^i$ exists, $1\leq i\leq k$.
\end{itemize}
\vskip 5mm
\begin{itemize}
\item $\mbox{\rm Re} (z)$: real part of $z\in \C$.
\item $\mbox{\rm Im}(z)$: imaginary part of $z\in\C$.
\end{itemize}

\end{document}